\documentclass[12pt,reqno]{amsart}
\usepackage[utf8]{inputenc}
\usepackage[english]{babel}
\usepackage{comment}
\usepackage{amsmath,amssymb,amsfonts,amsbsy,amsthm,amscd,latexsym,graphicx}
\usepackage{empheq,textcomp}
\usepackage{xcolor}
\usepackage{hyperref}

\theoremstyle{definition}
\newtheorem{theorem}{Theorem} [section]
\newtheorem{corollary}[theorem]{Corollary}
\newtheorem{lemma}[theorem]{Lemma}
\newtheorem{proposition}[theorem]{Proposition}
\newtheorem{definition}[theorem]{Definition}

\newtheorem{remark}[theorem]{Remark}
\newtheorem{example}[theorem]{Example}

\numberwithin{equation}{section}

\AtEndDocument{
  \par
  \bigskip
  \begin{tabular}{@{}l@{}}%
    \text{Department of Mathematics, University of Oregon, Eugene, OR 97403-1222, USA}\\
    \textit{E-mail address}: \texttt{putingyu@uoregon.edu}\\
  \end{tabular}}

\newcommand{\C}{\mathbb{C}}

\newcommand{\Fc}{{\mathcal{F}}}
\newcommand{\Gc}{{\mathcal{G}}}
\newcommand{\N}{\mathbb{N}}
\newcommand{\R}{\mathbb{R}}

\newcommand{\T}{\mathbb{T}}

\newcommand{\Z}{\mathbb{Z}}
\newcommand{\Q}{\mathbb{Q}}

\newcommand{\Rk}{{\mathfrak{R}}}

\newcommand{\Eq}{\, = \,}

\newcommand{\Le}{\, \le \,}

\newcommand{\Ge}{\, \geq \,}

\newcommand{\qeddef}{{\quad $\diamondsuit$}}

\newcommand{\bigabs}[1]{\bigl|\,#1\,\bigr|}
\newcommand{\Bigabs}[1]{\Bigl|\,#1\,\Bigr|}

\newcommand{\ip}[2]{\langle\,#1,#2\,\rangle}
\newcommand{\bigip}[2]{\bigl\langle \,#1, \, #2 \,\bigr\rangle}

\newcommand{\norm}[1]{\|\,#1\,\|}
\newcommand{\bignorm}[1]{\bigl\|\,#1\,\bigr\|}
\newcommand{\Bignorm}[1]{\Bigl\|\,#1\,\Bigr\|}

\newcommand{\bigparen}[1]{\bigl(\,#1\,\bigr)}
\newcommand{\Bigparen}[1]{\Bigl(\,#1\,\Bigr)}

\newcommand{\set}[1]{\{#1\}}
\newcommand{\bigset}[1]{\bigl\{#1\bigr\}}
\newcommand{\Bigset}[1]{\Bigl\{\,#1\,\Bigr\}}

\newcommand{\clspan}{{\overline{\text{span}}}}

\newcommand{\inN}{_{n\in\N}}
\newcommand{\sumli}{\sum_{n=1}^\infty}

\newcommand{\al}{\alpha}

\begin{document}
\title{Operations that are incompatible with certain
systems of translates in $L^2(\R)$.}
\author{Pu-Ting Yu}

\date{\today}

\pagestyle{plain}
\maketitle
\begin{abstract}
   We say that closed subspace $M$ of $L^2(\R)$ admits a \emph{complete set of semi-regular a-translates} if there exist some $a>0$, finitely many functions $g_1,\dots,g_N$, some subsets $J_1,\dots,J_N$ of $\Z$ and some finite subsets $\set{\al_{1j}}_{j=1}^{K_1},\dots,\set{\al_{Nj}}_{j=1}^{K_N}$ of $\R$ such that $$M=\clspan{\bigset{g_i(\cdot-ak), ~g_i(\cdot-\al_{ij})\,|\,k\in J_i,1\leq j\leq K_i\,}}_{i=1}^N.$$
    Here $\cdot$ denotes a generic variable. In the first half of this paper, we study whether the properties of being closed under modulation, dilation, reflection or Fourier transform is compatible with the existence of a complete set of semi-regular $a$-translates in closed subspaces of $L^2(\R)$. Specifically, we prove that a closed subspace of $L^2(\R)$ does not admit a complete set of semi-regular $a$-translates if it is closed under modulation or if it is closed under dilation with respect to a scaling factor $b$ satisfying $|b|>1.$ We also show that no infinite-dimensional closed subspace of $L^2(\R)$ can simultaneously be closed under Fourier transform and admit a complete set of semi-regular $a$-translates with  $a^2\in \Q$, whereas for any $a>0$, there do exist closed subspaces that are closed under reflection and admit a complete set of semi-regular $a$-translates.

    In the second half of this paper, we prove that a closed subspace of $L^2(\R)$ does not admit a frame formed by a system of translates if it contains a closed subspace that is closed under modulation and contains a nonzero function in $M^1(\R)$.  In addition, we present related results concerning the incompatibility between being closed under Fourier transform and the existence of  frames or Schauder bases of translates in closed subspaces of $L^2(\R)$. All results in this half can be extended to $L^2(\R^d)$ for any $d>1.$
\end{abstract}

\section{Introduction}
A \emph{system of translates} is a subset $S$ of $L^2(\R)$ that consists of arbitrarily countable translates by real numbers of finitely many functions in $L^2(\R)$. That is, $$S=\set{g_i(\cdot-\alpha_{ij})\,|\,1\leq i\leq N,\,1\leq j\leq N_i }$$ for some $N\in \N$, some $g_1,\dots,g_N\in L^2(\R)$,  and some $N_i\in \N\cup\set{\infty}$. 
Here $\cdot$ denotes a generic variable. The functions $g_1,\dots,g_N$ involved in a system of translates are usually called \emph{the generators} of such a system.
If $\set{\alpha_{ij}\,|\,1\leq j\leq N_i, 1\leq i\leq N}$ is contained in $a\Z$ for some $a>0$, then we say $S$ is a \emph{system of regular a-translates}.  Due to their simple structures, the question as to what kind of systems of translates can be utilized to approximate functions in $L^2(\R)$, in various senses, has been a major research topic in the fields of approximation theory, Banach space theory, operator theory, etc. For example, what kinds of systems of translates are complete in $L^2(\R)$? Equivalently, what kinds of systems of translates possess a closed span that is equal to $L^2(\R)$? The classic Fourier transform approach immediately leads us to the conclusion that no system of regular $a$-translates can be complete in $L^2(\R)$ whatever $a$ is. In the case that uncountably many translates of a function are considered, Wiener's Tauberian theorem (\cite{WI88}) states that the set of all translates of a function $g\in L^2(\R)$ is complete in $L^2(\R)$, i.e., $\clspan\bigset{g(\cdot-\alpha)\,|\,\alpha\in \R}$ is equal to $L^2(\R)$ if and only if the Fourier transform of $g$ is nonzero almost everywhere. The question regarding the completeness property of systems of translates becomes more intractable if we restrict the cardinality of associated translates to be countable, or even more, require the associated translates to be uniformly discrete. Here we say a countable subset $M\subseteq \R$ is \emph{uniformly discrete} if the distance between any two distinct numbers in $M$ is strictly bounded below by some positive number, i.e., $\inf_{x,y\in M}|x-y|>0.$ It was not known whether there exists a system of uniformly discrete translates that can be complete in $L^2(\R)$ until Olevskii constructed systems of uniformly discrete translates with a single generator that are complete in $L^2(\R)$ in \cite{OL96}. In particular, the associated translates $\set{\alpha_n}\inN$ taken in Olevskii's construction are nearly integers in the sense that $\alpha_n\notin \Z$ for all $n\in \N $ but $\lim\limits_{|n|\rightarrow \infty} |\alpha_n-n|\rightarrow 0.$ Later on, it was shown in \cite{OU04} that there exist complete systems of uniformly discrete translates generated by a single Schwartz function utilizing translates $\set{n+\alpha_n}\inN$ satisfying $0<|\alpha_n|<Cr^{|n|}$ for some $0<r<1$ and some $C>0.$ For relatively recent textbook recounting and other researches regarding the completeness of systems of translates in $L^2(\R)$, as well as in the setting of $L^p(\R)$ for $1\leq p<\infty$, we refer to \cite{AO96}, \cite{BOU06}, \cite{Le25}, \cite{LT25b}, \cite{OU16}  and \cite{OU18}.

Aside from the property of completeness, another related question along this direction is: under what conditions can a system of translates exhibit basis-like properties when appropriately ordered?
We say that a sequence $\set{x_n}\inN\subseteq L^2(\R)$ is a \emph{Schauder frame} for $L^2(\R)$ if there exists a sequence $\set{x_n^*}\subseteq L^2(\R)$, called the \emph{associated coefficient functionals}, such that 
\begin{equation}
    \label{basis_expansion}
x=\sumli\ip{x}{x_n^*}x_n
\end{equation}
with the convergence of the series in the norm of $L^2(\R)$ for all $x\in L^2(\R).$ Here the notation $\ip{\cdot}{\cdot}$ denotes the inner product associated with $L^2(\R)$. We call a Schauder frame an \emph{unconditional Schauder frame} for $L^2(\R)$ if the convergence of the series in Equation (\ref{basis_expansion}) is independent of order of the summation, i.e., $\sumli\ip{x}{x_{\sigma{(n)}}^*}\,x_{\sigma{(n)}}$ converges in the norm of $L^2(\R)$ for any permutation $\sigma:\N\rightarrow \N$. If the coefficient functionals associated with a Schauder frame (resp.\ unconditional Schauder frame) satisfies $\ip{x_m}{x_n^*}=\delta_{mn}$, then we call the Schauder frame a \emph{Schauder basis} (resp.\ \emph{unconditional (Schauder) basis}). It was conjectured by 
Olson and Zalik in \cite{OZ92} that no systems of translates formed by a single generator can be a Schauder basis for $L^2(\R)$. Although the Olson-Zalik conjecture itself remains open as of as of the time of writing, many related results to Olson-Zalik conjecture have been proved. For example, Olson and Zalik showed that no systems of translates formed by a single generator can be an unconditional basis for $L^2(\R)$. Christensen, Deng and Heil proved that no systems of translates can be a \emph{frame} (see Section \ref{frame_incom_complete} for the definition) for $L^2(\R)$ in \cite{CDH99}. Recently, Lev and Tselishchev confirmed in \cite{LT25} that there does not exist an unconditional frame formed by a system of translates. However, there do exist Schauder frames of translates for $L^2(\R)$(\cite{LT25}, \cite{FPT21}). The Olson-Zalik conjecture have been analogously formulated in the setting of $L^p(\R)$ for $1\leq p<\infty.$ For related research papers on this topic, we refer to \cite{FOSZ14}, \cite{LT23a}, \cite{LT24} and \cite{OSSZ11}.

While many results regarding systems of translates with certain approximation properties that $L^2(\R)$ can or cannot admit are known, there have been a few results concerning which types of closed subspaces of $L^2(\R)$ admit or fail to admit systems of translates with particular approximation properties. For example, it was shown in \cite{BDR94} that for a closed subspace $M$ of $L^2(\R)$ to admit a system of regular $1$-translates, it is necessary that there exist finitely many functions $g_1,\dots,g_N$ in $L^2(\R)$ such that 
\begin{equation}
\label{necessary_closed_sub_regular_trans}
\Fc(M)\subseteq \Bigset{\sum_{i=1}^N\sigma_ig_i\in L^2(\R)\,\big|\,\sigma_1,\dots,\sigma_N \text{ are 1-periodic functions defined on $\T$} }.
\end{equation}
Here $\Fc(M)$ denotes the set of Fourier transforms of all functions in $M$ and we identify $\T=\R/\Z$ with its fundamental domain, i.e., $\T=[0,1).$ Nevertheless, when an arbitrary closed subspace of $L^2(\R)$ is given, it is usually very hard to verify whether those $g_1,\dots,g_N$ in the inclusion relation (\ref{necessary_closed_sub_regular_trans}) do exist. In fact, those $g_1,\dots,g_N$ need not to be contained in $M$. We accordingly propose the following question: ``Are there any properties such that, if a closed subspace satisfies them, then it fails to admit a complete set of translates, an unconditional basis of translates, or a frame of translates, etc.?" Broadly speaking, we are seeking properties that are incompatible with the existence of systems of translates with certain approximation properties. As far as we are aware, not many researches have been conducted along this direction. Inspired by the results discovered in \cite{PY25}, this paper investigates whether the properties of being closed under modulation, dilation, reflection or Fourier transform (see Section \ref{preliminary} for the definitions of these operations) are compatible with the existence of complete sets of translates, unconditional bases of translates or frames of translates. Moreover, subspaces that satisfy these properties are of particular interests to us due to their primary roles in the fields of Gabor analysis, wavelet analysis, etc. Overall, our main results are summarized in the following two theorems.
\begin{theorem} \label{main_thm_1}
Let $M\subseteq L^2(\R)$ be an infinite-dimensional closed subspace. The following statements hold.
\begin{enumerate} \setlength\itemsep{0.5em}
    \item [\textup{(a)}] [Modulation] If there exists some real number $b\neq 0$ such that $$e^{2\pi ib\cdot} M \Eq\bigset{e^{2\pi ib\cdot}f(\cdot)\,|\,f\in M}\subseteq M,$$ 
    then $M$ does not admit a complete set of semi-regular $a$-translates for any $a>0.$

     \item [\textup{(b)}] [Dilation]  If there exists some nonzero real number $b$ with $|b|>1$ such that $$D_bM\Eq\bigset{f(bx)\,|\,f\in M}\subseteq M,$$
    then $M$ does not admit a complete set of semi-regular $a$-translates for any $a>0$. 

     \item [\textup{(c)}] [Fourier Transform] If $M$ is closed under Fourier transform, then $M$ does not admit a complete set of semi-regular $a$-translates for any $a>0$ satisfying $a^2\in \Q.$

\end{enumerate}
   Here a set of semi-regular $a$-translates is a set that has the form $$\bigcup_{i=1}^N\Bigparen{\bigset{g_i(\cdot-a k)\,\big|\,k\in J_i}\cup\bigset{g_i(\cdot-\alpha_ {ij})\,\big|\,1\leq j\leq K_i }}$$ 
   for some $N, K_1,\dots ,K_N\in \N$, some $g_1,\dots,g_N\in L^2(\R)$, some (possibly empty) subsets $J_1,\dots,J_N\subseteq \Z$ and some finite subsets $\set{\alpha_{1j}}_{j=1}^{K_1},\dots,\set{\alpha_{1j}}_{j=1}^{K_N}$ of $\R$.\qeddef
\end{theorem}

Consequently, numerous well-known closed subspaces in $L^2(\R)$ such as \emph {Hardy space} and those generated by \emph{Gabor systems} or \emph{wavelet systems} do not admit a complete set of semi-regular $a$-translates for any $a>0$ (See Section \ref{op_incom_complete} for examples).

As noted out above, it is known that $L^2(\R)$ does not admit any frame formed by a system of translates. However, to the best of our knowledge, there is no known characterization of closed subspaces of $L^2(\R)$ that fail to admit such frames. Our second main result provides such a characterization (See Section \ref{frame_incom_complete} for a brief introduction to $M^1(\R)$).
\begin{theorem}
\label{main_thm_2}
    Let $M\subseteq L^2(\R)$ be a closed subspace. Assume that $M$ contains a closed subspace $M'$ that is closed under modulation. That is, there exists some nonzero real number $b$ such that $\bigset{e^{2\pi ib\cdot}f(\cdot)\,|\,f\in M'}\subseteq M'$. Assume further that $M'$ contains a nonzero function in $M^1(\R)$.
    Then $M$ does not admit any frame formed by a system of translates (and hence, any unconditional basis formed by a system of translates).
\end{theorem}

This paper is organized as follows. Definitions, notations and some required known results are presented in Section \ref{preliminary}. Section \ref{op_incom_complete} is devoted to the proof of Theorem \ref{main_thm_1}. To be more precise, Section \ref{op_incom_complete} is divided into four parts. In each part, we investigate whether the existence of complete sets of semi-regular $a$-translates is compatible with one of the following four different operations on $L^2(\R)$: modulation, dilation, reflection and the Fourier transform. 
Finally, we prove Theorem \ref{main_thm_2} in Section \ref{frame_incom_complete}. The proof of Theorem \ref{main_thm_2} requires a variant of the classical Beurling density and some results related to modulation spaces, which will be briefly introduced in the beginning of Section \ref{frame_incom_complete}. Some partial results regarding the incompatibility between the existence of certain Schauder bases of translates and modulations are presented in this section as well.

\section{Preliminaries}
\label{preliminary}
 Throughout this paper, $L^2(\R)$ denotes the space of all Lebesgue square-integrable functions defined on $\R$ over the scalar field $\C.$
  Given any $a>0$, we identify $\frac{1}{a}\T=\R/a\Z$ with its fundamental domain, i.e.,  $\frac{1}{a}\T=[0,1/a).$ 
  We say that a complex-valued function $\tau$ defined on $\R$ is $a$-periodic function if $\tau(x+a)=\tau(x)$. So, every $a$-periodic is uniquely determined by its behavior on $[0,a).$ We define the Fourier transform and inverse Fourier transform by $$\Fc(f)(\xi)=\widehat{f}(\xi)=\int_{\R}f(x)e^{2\pi ix\xi}\,dx\text{  and  } \Fc^{-1}(f)(\xi)=\int_{\R}f(x)e^{2\pi ix\xi}\,dx,$$ respectively.
 \begin{definition} Let $b$ be a nonzero real number.
 \begin{enumerate} \setlength\itemsep{0.5em}
  \item [\textup{(a)}] The \emph{translation operator} $T_b\colon L^2(\R)\rightarrow L^2(\R)$ is the bounded linear operator defined by $(T_bf)(x)=f(x-b).$ 
     \item [\textup{(b)}] The \emph{modulation operator} $M_b\colon L^2(\R)\rightarrow L^2(\R)$ is the bounded linear operator defined by $(M_bf)(x)=e^{2\pi ibx}f(x).$ 
      \item [\textup{(c)}] For any nonzero real number $b$ with $b\neq 1$, the \emph{dilation operator} $D_b\colon L^2(\R)\rightarrow L^2(\R)$ is the bounded linear operator defined by $(D_bf)(x)=f(bx).$ 
       \item [\textup{(c)}] The \emph{reflection operator} $\Rk\colon L^2(\R)\rightarrow L^2(\R)$ is the bounded linear operator defined by $(\Rk f)(x)=f(-x).$ \qeddef
    \end{enumerate}   

 \end{definition}

 We adopt the convention that $$(M_bT_af)(x)=M_b\bigparen{(T_af)(x)}=e^{2\pi ibx}f(x-a)\Eq e^{2\pi iab}(T_aM_bf)(x),$$
 for all $a,b\in\R.$ The same convention applies to all other compositions of the operators defined above.
Throughout this paper, we say that a closed subspace $M\subseteq L^2(\R)$ is \emph{closed under modulation} if there exists some nonzero $b\in \R$ such that $M_bf\in M$ for all $f\in M$. Also, given a subset of $L^2(\R)$, we write $M_b(M)$ to mean the set $\set{e^{2\pi ib\cdot}f\,|\,f\in M}$. Same logic applies to other operations as well.

\begin{definition} Fix $a>0$ and let $M$ be a nontrival closed subspace of $L^2(\R)$, i.e., $M\neq \set{0}$.
\begin{enumerate} \setlength\itemsep{0.5em}

  \item [\textup{(a)}] \emph{A system (or set) of regular $a$-translates}, denoted by $S_{a,N}\bigparen{g_i,J_i}$, is a set that has the form 
  $$S_{a,N}(g_i,J_i)\Eq \set{T_{ak}g_i\,|\,k\in J_i,\,1\leq i\leq N}$$ for some $N\in \N$, some $g_1,\dots,g_N\in L^2(\R)$ and some $J_1,\dots,J_N\subseteq \Z.$ In the case that $J_1=\dots=J_N=J$, we will simply write $S_{a,N}(g_i,J)$. 

 \item [\textup{(b)}] \emph{A set of semi-regular $a$-translates}, denoted by $S_{a,N}\bigparen{g_i,J_i,\set{\alpha_{ij}}_{j=1}^{K_i}}$, is a set that has the form 
  $$S_{a,N}\bigparen{g_i,J_i,\set{\alpha_{ij}}_{j=1}^{K_i}}\Eq S_{a,N}\bigparen{g_i,J_i}\cup \set{T_{\alpha_{ij}}g_i\,|\,\,1\leq j\leq K_i,\,1\leq i\leq N}$$ for some $N,K_1,\dots ,K_N\in \N$, some $g_1,\dots,g_N\in L^2(\R)$, some $ J_1,\dots,J_N\subseteq \Z$ and some finite subsets $\set{\alpha_{1j}}_{j=1}^{K_1},\dots,\set{\alpha_{1j}}_{j=1}^{K_N}$ of $\R.$

  \item [\textup{(c)}] We say that $M$ admits a \emph{complete set of regular $a$-translates} (resp. complete set of semi-regular $a$-translates) if there exists a complete set of regular $a$-translates (resp. complete set of semi-regular $a$-translates) whose closed span equals to $M.$
    \end{enumerate}   
 In either of the definitions above, those functions $g_i$ are called \emph{generators} of the corresponding system. \qeddef
\end{definition}
 As mentioned in the Introduction, a necessary condition for a closed subspace $M$ to admit a complete set of regular $1$-translates is that the inclusion relation (\ref{necessary_closed_sub_regular_trans}) is satisfied. This necessary condition comes from the characterization of shift-invariant subspaces of $L^2(\R)$ using associated range functions, which we define as follows. 
 For more background knowledge on range functions, as well as characterizations of shift-invariant subspaces, we refer to \cite{He64}, \cite{BDR94} and \cite{Bo00}.
 
 \begin{definition} Fix $a>0$. 
\begin{enumerate}
 \item [\textup{(a)}] The \emph{fiberization map} (with respect to $a$) $\Psi_a:L^2(\R)\rightarrow L^2(a\T,\ell^2(\Z))$ is defined by $$\bigparen{\Psi_a(f)}(\xi)\Eq \bigparen{\widehat{f}(\xi-ak)}_{k\in \Z}$$
for almost every $\xi\in a\T.$ Here $L^2(a\T,\ell^2(\Z))$ is the Hilbert space of all $\ell^2(\Z)$-valued square-integrable functions defined on $a\T.$
    \item [\textup{(b)}] A \emph{range function} (of order $a$) is a mapping that maps $a\T$ to closed subspaces of $\ell^2(\Z).$ We say that a range function $J_a$ is (weakly) \emph{measurable} if for each $a,b\in \ell^2(\Z)$ the mapping $T\colon a\T\rightarrow \C$ defined by $T(x)=\ip{P(x)a}{b}$ is a measurable function, where $P(x)$ is the orthogonal projection from $\ell^2(\Z)$ onto $J_a(x).$  
    \item [\textup{(c)}] Let $J_a$ be a range function. The associated \emph{dimension function} with $J_a$ is the $a$-periodic function $V\colon a\T\rightarrow \N\cup\set{0,\infty} $ defined by $V(x)=\text{dim}(J_a(x)).$
\qeddef
\end{enumerate}
\end{definition}
 Clearly, for each $a>0$ the corresponding fiberization map is a bounded bijective linear operator from $L^2(\R)$ onto $L^2(a\T,\ell^2(\Z))$. We will see in the following proposition that by using the idea of fiberization maps along with the notion of range functions, we can obtain a characterization of all $a$-shift-invariant subspaces of $L^2(\R)$. Here we say that a closed subspace $M$ of $L^2(\R)$ is \emph{$a$-shift-invariant} if $T_{ak}(M)\subseteq M$ for all $k\in \Z$. We remark that the $1$-shift invariance version of the following proposition can be found in \cite{He64} (with slight modifications) and in \cite{Bo00}. Using the standard dilation technique, we can extend this result from the setting of $1$-shift invariance to $a$-shift invariance for any $a>0.$ For the sake of completeness, we include the short proof here. Nevertheless, we will primarily utilize the $L^2(\R)$ version of this result, we present it in that setting as follows.  
 \begin{proposition} \label{necessary_cond_range}
 Fix $a>0$.  The following statements hold.
\begin{enumerate}
    \item [\textup{(a)}]
    A closed subspace $M$ of $L^2(\R)$ if $a$-shift-invariant if and only if there exists a measurable range function $J_a$ for which $\bigparen{\Psi_{\frac{1}{a}}(f)}(\xi)\in J_a(\xi)$ for almost every $\xi\in \frac{1}{a}\T$ and all $f\in M.$ The correspondence between $a$-SI subspaces and range functions is one-to-one under the convention that two range functions are considered equal if they agree almost everywhere.
    \item [\textup{(b)}]The correspondence between $a$-shift-invariant subspaces and range functions is one-to-one under the convention that we consider two range functions are equal if they are equal almost everywhere.
    \item [\textup{(c)}] If there exists a countable subset $J$ of $\N$ and a countable subset $\set{g_i}_{i\in J}$ of $L^2(\R)$ such that $M=\clspan{\set{T_{ak}g_i\,|\,k\in \Z,\,i\in J}},$
     then the associated range function $J_a$ with $M$ is given by $J_a(\xi)=\clspan{\bigset{\bigparen{\Psi_{\frac{1}{a}}(g_i)}(\xi)\,|\,i\in J}}$.  Moreover, we have 
     $$M=\set{f\in L^2(\R)\,|\,\Psi_\frac{1}{a}(f)(\xi)\in J_a(\xi) \text{ for almost every $\xi\in \frac{1}{a}\T.$}}$$
    \end{enumerate}
    \begin{proof}
        Let $M$ be a $a$-shift-invariant subspace of $L^2(\R)$ and define the closed subspace $M'$ of $L^2(\R)$ by $M'\Eq \set{f(ax)\,|\,f\in M}.$ Take any $f\in M$ and let $h(x)=f(ax)$. Since $(T_{ak}f)(ax)=(T_kh)(x)$ for all $k\in \Z$, we see that $M$ is $a$-shift invariant if and only if $M'$ is $1$-shift-invariant. Let $J_1$ be the range function associated with $M'$ and define $J_a(\xi)=J_1(a\xi)$. Statements (a)--(c) then simply follow from \cite[Proposition 1.5]{Bo00}. 
        \end{proof}
 \end{proposition}
  An immediate corollary of Proposition \ref{necessary_cond_range} is the following necessary condition for a closed subspace $M$ of $L^2(\R)$ to admit a complete set of regular $a$-translates for some $a>0$.
 \begin{corollary}
     \label{nece_condit_complete_regu_trans} Fix $a>0.$ Assume that $M\subseteq L^2(\R)$ is a closed subspace that admits a complete set of regular $a$-translates $S_{a,N}(g_i,J_i)$. Then we have 
     $$\Fc(M)\subseteq \Bigset{f=\sum_{i=1}^N\sigma_i\widehat{g_i}\in L^2(\R)\,|\,\sigma_1,\dots,\sigma_N \text{ are }\frac{1}{a}\text{-periodic functions that are finite a.e}}.$$ \qeddef
 \end{corollary}
  Note that those $\frac{1}{a}$-periodic functions $\sigma_i$ associated with a given $f$ are not necessarily unique and $\sigma_i\widehat{g}_i$ needs not be in $L^2(\R).$ 
  
  Finally, we present a lemma which characterizes all generators of a given $a$-shift-invariant subspace with a single generator. The proof is a straightforward application of Corollary \ref{nece_condit_complete_regu_trans} and is essentially the proof of \cite[Corollary 2.4]{BDR94}. So, we skip the proof below.
  \begin{corollary}
  \label{generator_lemma}
      Fix $a>0$ and fix $g\in L^2(\R)$. Assume that $f\in S_{a,1}(g,\Z).$ Then $$S_{a,1}(f,\Z)=S_{a,1}(g,\Z) \text{  if and only if  }\text{supp}(\widehat{f})\Eq\text{supp}(\widehat{g}),$$
      where supp$(\widehat{f})$ denotes the support of $\widehat{f}$ in $\R$.
  \end{corollary}
 
\section{Operations That Are Incompatible With Complete Sets of Semi-regular $a$-translates}
\label{op_incom_complete}
In this section, we investigate whether the existence of a complete set of semi-regular $a$-translates in a closed subspace is compatible with the property of being closed under the following four operations: modulation, dilation, reflection and Fourier transform. We begin with the case of modulation, then move on to dilation, reflection, and then conclude this section on the case of the Fourier transform. 

\subsection{Modulation} Let $M\subseteq L^2(\R)$ be a nontrivial closed subspace of $L^2(\R)$ that is closed under modulation $M_b$ for some $b\neq0.$ We observe that if $M$ admits a complete set of regular $a$-translates $S_{a,N}(g_i,J_i)$, then for each $1\leq i\leq N$ we have 
$T_b(\widehat{g_i})=\sum_{k=1}^N\sigma_{ik}\widehat{g_k}$ 
for some $g_1,\dots,g_N$ and some $\frac{1}{a}$-periodic functions $\sigma_{i1},\dots,\sigma_{iN}$. So, the existence of a complete set of regular $a$-translates in $M$ depends on whether there exist $g_1,\dots,g_N$ satisfying the system of equations $T_b(\widehat{g_i})=\sum_{k=1}^N\sigma_{ik}\widehat{g_k}$ for $i=1,\dots,N$. We investigate the solvability of such equations first.
\bigskip

\noindent \textbf{Solvability of $a$-periodic Equations of Translates I}
\begin{definition}
    Fix $a>0$ amd fix $N\in \N$. Let $\sigma=\set{\sigma_i}_{i=1}^N\cup\set{\sigma_{ij}}_{1\leq i,j\leq N}$ be a collection of $a$-periodic functions that are finite almost everywhere and let $\Gamma=\set{\alpha_i}_{i=1}^N\cup \set{\alpha_{ij}}_{1\leq i,j\leq N}$ be a set of scalars. We say that the \emph{system of $a$-periodic equations of translates  with respect to $\sigma$ and $\Gamma$ is solvable} on $E$  if there exist some $g_1,\dots,g_N\in L^2(\R)$, not all identically zero, such that 
    \begin{equation}
        \label{a_periodic_eq_trans}
\sigma_{i}T_{\alpha_{i}}g_i=\sum_{j=1}^N\sigma_{ij}T_{\alpha_{ij}}g_j
        \end{equation}
    holds almost everywhere on $E+a\Z$ for all $1\leq i\leq N.$\qeddef
\end{definition}
Clearly, Equations (\ref{a_periodic_eq_trans}) are always satisfied if every $a$-periodic function in $\sigma$ is a zero function. Consequently, we will always implicitly assume that at least one of the equations in (\ref{a_periodic_eq_trans}) is nontrivial.
To the best of our knowledge, only a few results are known in this direction. The example presented below shows that for any $n_0,n_1\in \Z$ and any finite $a$-periodic function $\sigma_1$ the $a$-periodic equation of translates with respect to $\sigma=\set{1}\cup\set{\sigma_1}$ and $\Gamma=\set{an_0}\cup \set{an_1}$ is not solvable on $\R.$
\begin{example}
    Let $n_0,n_1\in \Z$ be some arbitrary integers and let $\sigma$ be a nonzero $a$-periodic function that is finite almost everywhere. Suppose that there exists some nonzero $g\in L^2(\R)$ such that 
    \begin{equation}
    \label{iterative_trans_regu_eq_order1}
        T_{an_0}g\Eq \sigma_1 T_{an_1}g
    \end{equation}
    Equation (\ref{iterative_trans_regu_eq_order1}) implies that $g$ is not identically zero on $[0,a].$
     Iterating Equation (\ref{iterative_trans_regu_eq_order1}), we see that for all $k\in \N$ $T_{ak(n_0-n_1)}g=\sigma_1^kg.$ Thus, we have \begin{align}
         \begin{split}
             \int_\R|g(x)|^2\,dx&\geq \sum_{k\in \N} \int_0^a|(T_{ak(n_0-n_1)}g)(x)|^2\,dx\\
             &\Eq \sum_{k\in \N} \int_0^a |\sigma_1(x)|^{2k}|g(x)|^2\,dx,
         \end{split}
     \end{align}
     which implies $|\sigma_1(x)|<1$ for almost every $x\in \R$. However, this would lead us to the conclusion that $\norm{T_{an_0}g}_{L^2(\R)}<\norm{T_{an_1}g}_{L^2(\R)}$, which is a contradiction. \qedhere
\end{example}
We will see that if we restrict all $\sigma_i$ in Equation (\ref{a_periodic_eq_trans}) to be equal to 1 almost everywhere and require all translates in Equation (\ref{a_periodic_eq_trans}) to lie in $a\Z$ for some $a>0$, then the corresponding $a$-periodic equations of translates is not solvable on $\R$. We need several lemmas to establish this result. 

We mention that the proof of the lemma below is the proof of \cite[Lemma 4.1]{BS10} with slight modifications. 
\begin{lemma}
\label{finite_generate_lemma} Fix $a>0$. 
    Let $M$ be a nontrivial $a$-shift-invariant subspace of $L^2(\R)$. Assume that $M=S_{a,N}(g_i,\Z)$ for some $g_1,\dots,g_N\in L^2(\R)$. 
    Then any nontrivial $a$-shift-invariant subspace of $M$ is not closed under modulation.
\end{lemma}
\begin{proof}
    Assume that $M'$ is a $a$-shift-invariant subspace of $M$ for which $M_b(M')\subseteq M'$ for some nonzero $b\in \R.$
    Let $J_{M'}(x)$ be the range function associated with $M'$. Then we define $\sigma_{M'}:\frac{1}{a}\T\rightarrow [0,1]$ by $$\sigma_{M'}(x+\frac{k}{a})=\norm{P_{M'}(x)e_k}_{\ell^2(\Z)}^2$$
for any $k\in \Z$ and $x\in \frac{1}{a}\T.$ Here $\set{e_k}_{k\in \Z}$ denotes the standard orthonormal basis for $\ell^2(\Z)$ and $P_{M'}(x)$ is the orthogonal projection from $\ell^2(\Z)$ onto $J_{M'}(x).$ Clearly, the dimension of $J_{M'}(x)$ is equal to $\sum_{k\in \Z}\sigma_{M'}(x+\frac{k}{a})$ for almost every $x\in \frac{1}{a}\T.$ Since $M_b(M')\subseteq M'$, we have $\sigma_{M'}(x+b)\geq \sigma_{M'}(x)$ for almost every $x\in \R$ by \cite[Proposition 2.6]{BR03}. Let $E$ be the support of $\sigma_{M'}$. By iterating the inequality $\sigma_{M'}(x+b)\geq \sigma_{M'}(x)$, we obtain that $\sum_{k\in \Z}\sigma_{M'}(x+bk)=\infty$ on $E+b\Z$. Here $E+b\Z$ means $\cup_{k\in\Z} (E+bk)$. Let $L$ be the minimum number of generators of $M$. Then we compute 
    \begin{align}
    \begin{split}
    \frac{L}{a}\Ge \int_0^{1/a} \dim (J_{M'}(x))\,dx &\Eq \int_0^{1/a} \sum_{k\in \Z}\sigma_{M'}(x+\frac{k}{a})\,dx \\&\Eq \int_{\R}\sigma_{M'}(x)\,dx\\
    &\Eq \int_0^b \sum_{k\in \Z}\sigma_{M'}(x+bk)\,dx\\
    &\Eq \infty,
    \end{split}
\end{align}
which is a contradiction.
\end{proof}
\begin{remark}
    The function $\sigma_{M'}$ defined in the proof of Lemma \ref{finite_generate_lemma} is called the \emph{spectral function} associated with $J_{M'}$, which was first introduced and studied by Rzeszotnik in \cite{RZ00}. For more results regarding spectral functions, we refer to \cite{BR03}.
\end{remark}

The next lemma and its proof comes from the key idea of the proof of \cite[Theorem 3.4]{BS10} and \cite[Theorem 4.2]{BS10}.
\begin{lemma}\label{linear_independence_lattice}
Fix $a,b>0$ and fix $N\in \N.$ Let $\set{g_i}_{i=1}^N \subseteq L^2(\R)$ be a finite collection of nonzero functions and let $J_1,\dots,J_N$ be some subsets of $\Z.$ Assume that the closed subspace $$M= \clspan{\set{M_{bn}T_{ak}g_i\,|\,n\in J_i,\,k\in \Z}}_{i=1}^N$$ admits a complete set of regular $a$-translates. Then the following statements hold.
\begin{enumerate} \setlength\itemsep{0.5em}
    \item [\textup{(a)}] If $\cup_{i=1}^NJ_i$ is bounded below, then $M_{bn_0}S_{a,N}\bigparen{g_i,\Z}\nsubseteq M$ for any $n_0< \inf(\cup_{i=1}^NJ_i).$  
     \item [\textup{(b)}] If $\cup_{i=1}^NJ_i$ is bounded above, then $M_{bn_0}S_{a,N}\bigparen{g_i,\Z}\nsubseteq M$ for any $n_0> \sup(\cup_{i=1}^NJ_i).$ 
\end{enumerate}
\end{lemma}
\begin{proof}
By applying the reflection operator if necessary, we may assume that $J$ is bounded below. Take any  $n_0<\inf\set{n\in  J_i\,|\,1\leq i\leq N}$ and suppose that $M_{bn_0}g_i\in  M$ for all $1\leq i\leq N.$ For any $k_1\leq k_2\in \Z\cup \set{-\infty,\infty}$ we define the $a$-shift invariant subspace $$V_{k_1,k_2}=\clspan{\set{M_{bj}T_{ak}g_i\,|\,k_1\leq j\leq k_2,\,\,k\in \Z}}_{i=1}^N.$$
By $V_{-\infty,k_1}$ and $V_{k_1,\infty}$ we mean $\displaystyle\bigcup_{k\leq k_1}V_{k,k_1}$ and $\displaystyle\bigcup_{k\geq k_1}V_{k_1,k}$, respectively. Let $L=\sup(\cup_{i=1}^NJ_i)$. Note that $L$ might be $\infty.$
    By assumption, we have $M_{bn_0}g_i\in V_{n_0,L}\subseteq V_{n_0+1,L}.$ Consequently, we have $V_{n_0,L}=V_{n_0+1,L}$. Applying the modulation operator $M_{bk}$ to the equation $V_{n_0,L}=V_{n_0+1,L}$ for each $k\in \Z$, it follows that 
    \begin{equation}
     \label{nesting_equation}
    V_{n_0+k,L+k}=V_{n_0+1+k,L+k},
    \end{equation}
    for any $k\in \Z.$ Next, we show that $V_{n_0-r,L}\subseteq V_{n_0+1,L}$ for any integer $r\geq 0$ by using mathematical induction. We have proved the base case $r=0$. Assume that $V_{n_0-r,L}\subseteq V_{n_0+1,L}$ for some $r\geq 0.$ Since $V_{n_0-r-1,\,L-r-1}\subseteq V_{n_0-r,\,L-r-1}$ by Equation (\ref{nesting_equation}), we see that $$V_{n_0-r-1,\,L}\subseteq V_{n_0-r-1,\,L-r-1}\cup V_{n_0-r,L}\subseteq V_{n_0-r,\,L-r-1} \cup V_{n_0+1,L}\subseteq V_{n_0+1,L}.$$
    By mathematical induction, we see that  $V_{n_0-r,L}\subseteq V_{n_0+1,L}$ for any $r\geq 0$. Therefore, we have $$V_{-\infty,L}\Eq  \bigcup_{r=0}^\infty V_{n_0-r,L}\Eq  V_{n_0+1,L},$$ which implies that $V_{-\infty,L}$ is contained in a $a$-shift-invariant closed subspace generated by finitely many functions. 
    However, since $$V_{-\infty,L}\subseteq V_{-\infty,L+1}= M_{b}(V_{-\infty,L}),$$
    we obtain a contradiction to Lemma \ref{finite_generate_lemma}.
\end{proof}
\begin{theorem}
\label{Non_exist_iterative_eq}
Fix $a,b>0$ and fix $N\in \N$.  Let $\sigma=\set{1}_{i=1}^N\cup\set{\sigma_{ij}}_{1\leq i,j\leq N}$ be a collection of $b$-periodic functions that are finite almost everywhere and let $\Gamma=\set{an_0}_{i=1}^N\cup\set{an_{ij}}_{1\leq i,j\leq N}$ be a set of scalars with $n_0\notin [\min\limits_{n\in \Gamma} n,\,\max\limits_{n\in \Gamma} n ]$. Then the $b$-periodic equations of translates with respect to $\sigma,\Gamma$ are not solvable on $\R$. 

That is, for any $g_1,\dots,g_N\in L^2(\R)$, not all identically zero, there exists at least one $1\leq i\leq N$ for which $$T_{an_0}g_{i}\neq \sum_{j=1}^{N}\sigma_{ij}T_{an_{ij}}g_j$$ 
on a set $E\subseteq \R$ of positive measure.
\end{theorem}
\begin{proof} 
Suppose there exists $g_1,\dots,g_N$ in $L^2(\R)$, not all identically zero, such that 
\begin{equation}
       \label{iterative_equation}
      T_{an_0}g_i\Eq \sum_{j=1}^{N}\sigma_{ij}T_{an_{ij}}g_{j}, 
    \end{equation}
    on $\R$ for all $1\leq i\leq N$. Let $m(x)= \max\limits_{1\leq i,j\leq N}\sigma_{ij}(x)$. Then we define a $b$-periodic function $\tau$ on $[0,b]$ by \begin{equation*}
\begin{aligned}
\tau(x)=\left\{
             \begin{array}{ll}
             \frac{1}{m(x)}, &\text{if $m(x)>1$} \\
             
            1, &\text{if $m(x)\leq 1$.}
             \end{array}
\right.
\end{aligned}
\end{equation*}
Then $\tau(x)$ is nonzero almost everywhere. Moreover, $\tau\sigma_{ij}g_j$ is in $L^2(\R)$ for all $1\leq i,j\leq N$.  Fix $1\leq \ell\leq N.$ Multiplying Equation (\ref{iterative_equation}) by $\tau$, we obtain  
\begin{equation}
       \label{iterative_equation_L2}
      \tau T_{an_0}g_\ell\Eq \sum_{j=1}^{N}\tau_{ij}T_{an_{ij}}g_{j}, 
    \end{equation}
    where $\tau_{ij}=\tau\sigma_{ij}.$ Taking the inverse Fourier transform of Equation (\ref{iterative_equation_L2}), we rewrite Equation (\ref{iterative_equation_L2}) as
    \begin{equation}
       \label{iterative_equation_L2_2}
      M_{an_0}\Fc^{-1}\bigparen{(T_{-an_0}\tau)g_{\ell}}\Eq \sum_{j=1}^{N}M_{an_{ij}}\Fc^{-1}\bigparen{(T_{-an_{ij}}\tau_{ij})g_{j}}.
    \end{equation}
Note that $\Fc^{-1}\bigparen{(T_{-an_{ij}}\tau_{ij})g_{j}}$ is contained in $\clspan{\bigset{S_{\frac{1}{b},1}(\Fc^{-1}(g_j),\Z)}}$ by Corollary \ref{nece_condit_complete_regu_trans}. So, Equation (\ref{iterative_equation_L2_2}) implies that 
    \begin{equation}
       \label{iterative_equation_L2_3}
    M_{an_0}S_{\frac{1}{b},1}\bigparen{\Fc^{-1}\bigparen{(T_{-an_0}\tau) g_{\ell}},\Z}\subseteq \clspan{\bigset{M_{an_{ij}}T_{bk}(\Fc^{-1}g_j)\,|\,1\leq i,j \leq N,\,k\in\Z }}.
    \end{equation}
    On the other hand, since $\text{supp}\bigparen{(T_{-an_0}\tau)g_\ell}=\text{supp}(g_\ell)$, we also have, by Corollary \ref{generator_lemma}, $$\clspan{\bigset{S_{\frac{1}{b},1}\bigparen{\Fc^{-1}g_{\ell},\Z}}}\Eq \clspan{\bigset{S_{\frac{1}{b},1}\bigparen{\Fc^{-1}\bigparen{(T_{-an_0}\tau) g_{\ell}},\Z}}} $$
    Consequently, we obtain that  
    $$M_{an_0}S_{\frac{1}{b},1}\bigparen{\Fc^{-1} g_{\ell},\Z}\subseteq \clspan{\bigset{M_{an_{ij}}T_{bk}(\Fc^{-1}g_j)\,|\,1\leq i,j \leq N,\,k\in\Z }}. $$
   Since $\ell$ is arbitrary, our assumption then leads us to the conclusion that 
    $$\ M_{an_0}S_{\frac{1}{b},N}\bigparen{\Fc^{-1} g_{\ell},\Z}\subseteq \clspan{\bigset{M_{an_{ij}}T_{bk}(\Fc^{-1}g_j)\,|\,1\leq i,j \leq N,\,k\in\Z }}, $$
    which is a contradiction to Lemma \ref{linear_independence_lattice}.
\end{proof}
\begin{remark}
The \emph{HRT (Heil-Ramanathan-Topiwala) Conjecture} \cite{HRT96} states that any finite set of time-frequency shifts of a non-zero function $g\in L^2(\R)$ is linearly independent. That is, for any $N\in \N$, if 
\begin{equation}
\sum_{n=1}^Nc_nM_{w_n}T_{\gamma_n}g=0
\end{equation}
for some finite collection $\set{(\gamma_n,w_n)}_{n=1}^N\subseteq \R^2$, then we must have $c_n=0$ for all $1\leq n\leq N.$ This conjecture remains open at the time of writing. A known results regarding the HRT conjecture due to Linnell \cite{Li99} shows that the HRT conjecture holds if $\set{(\gamma_n,w_n)}_{n=1}^N$ lies in $a\Z\times b\Z$ for some $a,b>0.$  Using the language of $a$-periodic equations of translates, Linnell's result is equivalent to the statement that for any increasing sequence of integers $\set{m_i}_{i=0}^N$ and any sequence of integers $\set{n_i}_{i=0}^N$ the following $a$-periodic equations of translates of order $N$, $$M_{an_0}T_{bm_0}g=\sum_{j=1}^N c_jM_{an_j}T_{bm_j}g,$$
is not solvable on $\R$ for any scalars $c_1,\cdots,c_N$. So, Theorem \ref{Non_exist_iterative_eq} not only extends the lattice case of HRT conjecture to a more general setting of generic periodic functions that are finite a.e., but also provides insight what form HRT conjecture might be formulated when considering the linear independence of time-frequency shifts of more than one function. \qeddef 
\end{remark}
\bigskip

\noindent \textbf{Incompatibility Between Modulation and Complete Sets of Semi-regular Translates} We are now ready to prove that every closed subspace that is closed under modulation fails to admit a complete set of semi-regular $a$-translates for any $a>0.$ Note that every closed subspace of $L^2(\R)$ that is closed under modulation is infinite-dimensional. 
\begin{corollary}
\label{non_regu_trans_closed_modu}
 Let $M$ be a nontrivial closed subspace of $L^2(\R)$ that is closed under modulation. Then $M$ does not admits a complete set of regular $a$-translates for any $a>0.$
\end{corollary}
\begin{proof}
    Suppose to the contrary that there exist some $a>0$, $N\in \N$ and some functions $g_1,\dots,g_N$ in $M$ such that   $$M\Eq \clspan{\set{S_{a,N}(g_{i},J_i)}}\subseteq\clspan{\set{S_{a,N}(g_{i},\Z)}}.$$ 
    By considering $g_i=T_{ak_0}g_i$ for some $k_0\in\Z$ if necessary, we may assume that $g_i\in M.$
    Let $b$ be a nonzero real number for which $M_b(M)\subseteq M.$
    By Lemma \ref{nece_condit_complete_regu_trans}, for each $1\leq i\leq N$ there exist $\frac{1}{a}$-periodic functions $\sigma_{i1},\dots,\sigma_{iN}$ that are finite almost everywhere such that
    we have $$\Fc(M_bg_i)\Eq T_b\widehat{g_i}\Eq\sum_{j=1}^N\sigma_{ij}\widehat{g}_j,$$
    which is a contradiction to Theorem \ref{Non_exist_iterative_eq}. 
\end{proof}

\begin{theorem}
\label{modu_incomp_semiregu_trans}
      Let $M$ be a nontrivial closed subspace of $L^2(\R)$. Assume that $M$ is closed under modulation. Then $M$ does not admit any complete set of semi-regular $a$-translates. 
\end{theorem}
\begin{proof}
    
    Suppose to the contrary that there exist some $a>0$, some $N\in \N$, some nonzero functions $g_1,\dots,g_N$ in $L^2(\R)$, some subsets $J_1,\dots,J_N$ of $\Z$ and some finite subsets $\set{\alpha_{1j}}_{j=1}^{K_1},\dots,$ $\set{\alpha_{Nj}}_{j=1}^{K_N}$ of $\R$ such that $M=\clspan{\bigset{\bigcup_{i=1}^N S_a(g_i,J_i,\set{\alpha_{ij}}_{j=1}^{K_i})}}.$ Likewise, we may assume that $g_i\in M$ for all $1\leq i\leq N.$ Moreover, by Corollary \ref{non_regu_trans_closed_modu}, it suffices to consider the case that $V=\clspan{\set{T_{\alpha_{ij}}g_i\,|\,1\leq j\leq K_i,\,1\leq i\leq N}}\not\subseteq \clspan{\set{S_{a,N}(g_i,I_j)}}$ for any subsets $I_1,\dots,I_N$ of $\Z.$ Since $V$ is finite-dimensional, there exists some $\ell_1\in \N$ such that for any $\ell\geq \ell_1$ we have $M_{b\ell}g_j\notin V$ for all $1\leq j\leq N.$ Therefore, for each $1\leq j\leq N$ and each $\ell\geq \ell_1$ there exist some nonzero element $v_{j,\ell}\in V$ and some $\frac{1}{a}$-periodic functions $\sigma^{(\ell)}_{1j},\dots,\sigma^{(\ell)}_{Nj}$, not all identically zero, such that 
    $$T_{b\ell}\widehat{g_j}=\sum_{i=1}^N\sigma^{(\ell)}_{ij}\widehat{g_i}+v_{j,\ell}.$$
   Consequently, the closed subspace $\clspan{\bigset{T_{b\ell}\widehat{g_j}-\sum_{i=1}^N\sigma^{(\ell)}_{ij}\widehat{g_i}}_{\ell\geq \ell_1}}$ must be finite-dimensional for all $1\leq j\leq N.$ Therefore, there exists some $\ell_2\in \N$ such that for any $n\geq \ell_2+1$ we have 
    $$T_{bn}\widehat{g_j}- \sum_{i=1}^N\sigma^{(n)}_{ij}\widehat{g_i}\in \clspan{\bigset{T_{b{\ell}}\widehat{g_j}-\sum_{i=1}^N\sigma^{(\ell)}_{ij}\widehat{g_i}}^{\ell_2}_{\ell= \ell_1}}$$
    for all $1\leq j\leq N.$
    Hence, for each $1\leq j\leq N$ and each $n\geq \ell_2+1$ there exist some scalars $(d^{(n)}_{\ell,j})_{\ell_1\leq \ell\leq \ell_2}$ such that 
$T_{bn}\widehat{g_j}-\sum_{i=1}^N\sigma^{(n)}_{ij}\widehat{g_i}=\sum_{\ell=\ell_1}^{\ell_2}d^{(n)}_{\ell,j}\bigparen{T_{b\ell}\widehat{g_j}-\sum_{i=1}^N\sigma^{(\ell)}_{ij}\widehat{g_i}}.$ It follows that $$M_{bn}g_j\in \clspan{\set{M_{b\ell}T_{ak}g_i\,|\,\ell_1\leq \ell\leq \ell_2,\,k\in \Z,\,1\leq i\leq N}}$$ for all $1\leq j\leq N$. Equivalently, we have
    $$g_j\in \clspan{\set{M_{b(\ell-n)}T_{ak}g_i\,|\,\ell_1\leq \ell\leq \ell_2,\,k\in \Z,\,1\leq i\leq N}},$$
    for all $1\leq j\leq N$. Since $n\geq \ell_2+1$, we obtain a contradiction to Lemma \ref{linear_independence_lattice}.
\end{proof}
\begin{example}
(a) Let $E\subseteq \R$ be an arbitrary subset of positive measure. By Theorem \ref{modu_incomp_semiregu_trans}, the closed subspace $M=\set{f\in L^2(\R)\,|\,f=0~\text{on } E}$ does not admit any complete set of semi-regular $a$-translates. Equivalently, if $M=\set{f\in L^2(\R)\,|\,\Fc(f)=0 ~\text{on } E}$, then we have
$$M\neq \clspan{\bigset{\set{M_{bk_i}g_i\,|\,k_i\in J_i,~1\leq i\leq N}\cup \set{M_{\alpha_{ij}}g_i\,|\,1\leq i\leq N,~1\leq j\leq N_i}}},$$
for any $N,K_1,\dots ,K_N\in \N$, some $g_1,\dots,g_N\in L^2(\R)$, any $J_1,\dots,J_N\subseteq \Z$ and any finite subsets $\set{\alpha_{1j}}_{j=1}^{K_1},\dots,\set{\alpha_{1j}}_{j=1}^{K_N}$ of $\R$.

\smallskip
(b) Another example is the \emph{Weyl-Heisenberg systems} or \emph{Gabor systems} (see Section \ref{frame_incom_complete} for more introduction). By Theorem \ref{modu_incomp_semiregu_trans}, the closed span of any non-trivial Gabor system never admits a complete set of semi-regular $a$-translates for any $a>0$. For example, fix $\alpha,\beta>0$
and let $\phi=e^{-t}\chi_{[0,\infty)}$. It is known that the closed subspace 
$M=\clspan{\set{M_{\beta k}T_{\alpha n}\phi\,|\,n,k\in \Z}}$ is exactly $L^2(\R)$ if $\alpha\beta\leq 1$ and is a proper subspace of $L^2(\R)$ if $\alpha\beta>1$. By Theorem \ref{modu_incomp_semiregu_trans}, $M$ does not admit a complete set of semi-regular $a$-translates for any $a>0$ whatever $\alpha\beta$ is. \qeddef 
\end{example}
\subsection{Dilation}
\label{dilation_section} Next, we consider the case in which the closed subspaces are closed under dilation. 
The following theorem was pointed out to the author by Ryszard Szwarc, which shows that for any nonzero function $g\in L^2(\R)$ the set $\set{(D_{b})^ng}_{n\in \N\cup\set{0}}$ is finitely linearly independent for any $b>0$ with $b\neq 1$. By considering $b^2$ if necessary, we see that every closed subspace of $L^2(\R)$ that is closed under dilation is infinite-dimensional when $b\notin \set{0,1}.$ 
\begin{theorem}
\label{dila_fini_linear_inde}
    Fix $b>0$ with $b\neq 1$. Then for any $g\in L^2(\R)\setminus \set{0}$ the set $\set{D_{b^n}g}_{n\in \N\cup\set{0}}$ is finitely linearly independent in $L^2(\R)$.
\end{theorem}
\begin{proof}Let $g\in L^2(\R)$ be an arbitrary nonzero function. Suppose the contrary that there exist a nonzero sequence of scalars $(c_i)_{i=1}^N$ with $c_N= 1$ for which $$0=\sum_{i=1}^Nc_ig(b^ix).$$
This is equivalent to assume that the bounded linear operator $T\colon L^2(\R)\rightarrow L^2(\R)$ defined by $$T(f)=\Bigparen{\sum_{i=1}^Nc_iD_{b^i}}(f)$$ has a non-trivial kernel. Here the kernel of $T$ is the set $\set{f\in L^2(\R)\,|\,T(f)=0}.$
The operator $T$ can be factorized into a composition of $n$ linear factor of the form $(D_b-\lambda I)$ for some $\lambda\in \C$. Consequently, it suffices to show that for any $\lambda\in \C$ the bounded linear operator $D_b-\lambda I$ has a trivial kernel to reach a contradiction. Let $f\in L^2(\R)$ be a nonzero function such that $D_b(f)=\lambda f.$ A straightforward computation shows that $\lambda=b^{-1/2}$, and hence\begin{equation}
\label{dila_eq}
    f=b^{1/2}D_bf
\end{equation}
Iterating Equation (\ref{dila_eq}), we see that $f(x)=b^{k/2}f(b^{k}x)$ and $f(b^{-k}x)=b^{k/2}f(x)$ for all $k\in \N.$ Let $[c,d]$ be an arbitrary interval that is either contained in $(0,\infty)$ or $(-\infty,0)$.
Then we compute 
\begin{align}
    \begin{split}
        b^{k}\int_c^d|f(x)|^2\,dx&\Eq  \int_c^d|f(b^{-k}x)|^2\,dx\\
        &\Eq b^{k}\int_{b^kc}^{b^kd}|f(x)|^2\,dx.
    \end{split}
\end{align}
Since $\int_{b^kc}^{b^kd}|f(x)|^2\,dx\rightarrow 0$ as $k\rightarrow \infty,$ we see that $f$ must be the zero function. Therefore, the kernel of $D_b-\lambda I$ only contains the zero function.
\end{proof}

Next, we present examples that for each $a>0$ there exist closed subspaces of $L^2(\R)$ that admit a complete set of regular $a$-translates and are closed under dilation with the scaling factor $b$ satisfying $0<|b|\leq 1$. Using more than one generator, we can also show that there are closed subspaces of $L^2(\R)$ that
are closed under dilation with the scaling factor $b$ being irrational while still admitting a complete set of regular a-translates.

\begin{example}
\label{Dilation_compa_translates} Fix $a>0.$
    Let $g=\chi_{[0,a)}$ be the characteristic function defined on $[0,a).$ Then clearly, we have $$\bigset{\chi_{[am,an)}\,|\,n>m \text{ and } n,m\in \Z} \subseteq \clspan{\set{T_{ak}(\chi_{[0,a)})}_{k\in \Z}}.$$ 
    For each $k\in \Z$ and any nonzero real number $b$ with $\frac{1}{b}\in \Z$, we have $$D_b(T_{ak}(\chi_{[0,a)}))=D_b(\chi_{[ka,(k+1)a)})\in \bigset{\chi_{[ma,na)}\,|\,n>m \text{ and } n,m\in\Z}.$$
    Consequently, $\clspan{\set{T_{ak}(\chi_{[0,a)})}_{k\in \Z}}$ is closed under dilation $D_b$ for any nonzero real number $b$ with $\frac{1}{b}\in \Z.$ Even more, $\set{T_{ak}(\chi_{[0,a)})}_{k\in \Z}$ is an orthonormal basis for its closed span.\qeddef
\end{example}

\begin{example}
\label{Dilation_compa_translatesII} Let $b>0$ be such that $0<|b|\leq1$. Assume that there exists some integer $\ell\geq 2$ such that $Lb^\ell=b$ for some $L\in\Z.$ Clearly, we have $$\bigset{\chi_{[mb^{-j},nb^{-j})}\,|\,n>m \text{ and } n,m\in \Z} \subseteq \clspan{\set{T_{b^{-j}k}(\chi_{[0,b^{-j})})}_{k\in \Z}}.$$ 
    For each $k\in \Z$ and each $1\leq j\leq \ell-1$ we have $$D_b(T_{b^{-j}k}(\chi_{[0,b^{-j})}))\Eq D_b(\chi_{[kb^{-j},(k+1)b^{-j})})\Eq \chi_{[kb^{-j-1},(k+1)b^{-j-1})}\in M$$
Thus, the closed subspace $$M=\clspan{\set{T_{b^{-j}k}(\chi_{[0,b^{-j})})\,|\,1\leq j\leq \ell-1,\,k\in\Z}}.$$
is closed under the dilation operator $D_b.$ 
\qeddef
\end{example}

However, the result dramatically changes when the magnitude of the scaling factor exceeds $1$. We will show that no closed subspace of $L^2(\R)$ can be closed under dilation with the magnitude of the scaling factor being greater than $1$ while still admitting a complete set of semi-regular $a$-translates. We begin by establishing the following two lemmas.
\begin{lemma}\label{intersection_positive_measure}
 Fix $a>0,b>1$ and fix $N\in \N$. For any measurable subset $E\subseteq \R$ of positive measure there exists some $L(a,N)\in\N$  such that $|\bigcap_{j=0}^N(b^\ell E+a^{-1}j)|>0$ for all $\ell\geq L.$  
\end{lemma}
\begin{proof}
    Let $\epsilon= \frac{1}{2(N+1)}.$ By Lebesgue Density Theorem, there exists some interval $I$ such that $$|E\cap I|\geq (1-\epsilon)|I|.$$
    We then choose $L$ large enough that $b^{L}|I|\epsilon>a^{-1}N$ and let $\ell\geq L.$  Since $(b^\ell E\cap b^\ell I)\subseteq b^\ell E$, it suffices to prove that $|\bigcap_{j=0}^N\bigparen{(b^\ell E\cap b^\ell I)+a^{-1}j}|>0.$ 
    Observe that  
    \begin{align*}
        \begin{split}
 \Bigabs{\bigcap_{j=0}^N\bigparen{(b^\ell E\cap b^\ell I)+a^{-1}j}}&\Eq \Bigabs{\bigcap_{j=0}^N\bigparen{(b^\ell E+a^{-1}j)\cap (b^\ell I+a^{-1}j)}} \\
 &\Eq \Bigabs{\bigcap_{j=0}^N (b^\ell I+a^{-1}j)}-\Bigabs{\bigparen{\bigcap_{j=0}^N (b^\ell I+a^{-1}j)}\setminus\bigparen{\bigcap_{j=0}^N (b^\ell E+a^{-1}j)}}.
 \end{split}
 \end{align*}
We then compute 
\begin{align}
    \begin{split}
      \Bigabs{\bigparen{\bigcap_{j=0}^N (b^\ell I+a^{-1}j)}\setminus\bigparen{\bigcap_{j=0}^N (b^\ell E+a^{-1}j)}}&\Le\Bigabs{\bigcup_{j=0}^N\bigparen{(b^\ell I+a^{-1}j)\cap (b^\ell E+a^{-1}j)^c}}\\  
      &\Le \sum_{j=0}^N |(b^\ell I\setminus b^\ell E)+a^{-1}j|\\
      &\Eq \sum_{j=0}^N |b^\ell I\setminus b^\ell E|\\
      &\Le b^\ell(N+1) \epsilon |I|\\
      & \Eq \frac{b^\ell}{2}|I|.
    \end{split}
\end{align}
Since $\bigabs{\bigcap_{j=0}^N (b^\ell I+a^{-1}j)}= b^\ell|I|-Na^{-1}> (1-\epsilon)|I|$, it follows that 
$$\Bigabs{\bigcap_{j=0}^N\bigparen{(b^\ell E\cap b^\ell I)+a^{-1}j}}>b^\ell(1-\epsilon-\frac{1}{2})|I|>0.$$
\end{proof}

Given any $a>0$, recall that the fiberization map $\Psi_a\colon L^2(\R)\rightarrow L^2(a\T,\ell^2(\Z))$ is defined by $\bigparen{\Psi_a(f)}(\xi)\Eq \bigparen{\widehat{f}(\xi-ak)}_{k\in \Z}$. Hence, for any $b>0$, any $\ell\in\N$ and any $n_k\in \Z$ we have $$\Psi_{\frac{1}{a}}(T_{an_k}D_{b^{\ell}}f)(\xi)=\Bigparen{\widehat{f}\bigparen{b^{-\ell}(\xi-\frac{n}{a})}e^{-2\pi i an_kb^{-\ell}(\xi-\frac{n}{a})}}_{n\in \Z}$$
for all $f\in L^2(\R).$

\begin{lemma}
\label{dila_infini_dim_lemma}
    Fix $a>0, b>1$ and let $g\in L^2(\R)$ be a nonzero function. Assume that $J=\set{n_k}_{k\in \N}$ is an infinite subset of $\Z$. Then for each $N\in \N$ there exists $\ell(N,a)\in \N$ and some subset $I\subseteq \frac{1}{a}\T$ of positive measure such that $$\text{dim}\bigparen{\clspan\set{\Psi_{\frac{1}{a}}(T_{an_k}D_{b^\ell}g)(\xi)}}_{k=1}^N\Eq N$$
    for almost every $\xi\in I.$
\end{lemma}
\begin{proof} Let $E$ be a measurable subset of positive measure that is contained in the support of $\widehat{g}.$ The case $N=1$ clearly holds. Now fix $N\geq 2$ and pick $\ell$ large enough so that the following statements hold:
\begin{enumerate} \setlength\itemsep{0.5em}
    \item [\textup{(a)}] $\{e^{2\pi i(n_k-n_N)b^{-\ell}}\}_{k=1}^{N-1}$ are $N-1$ distinct complex numbers,

     \item [\textup{(b)}] $e^{2\pi i(n_k-n_N)b^{-\ell}}\neq 1$ for all $1\leq k\leq N-1,$
\item [\textup{(c)}] $|\bigcap_{m=0}^{N-1}(b^\ell E+a^{-1}m)|>0$. (By Lemma \ref{intersection_positive_measure})    
\end{enumerate}
Note that statement (c) implies that  $$\Fc(D_{b^\ell}g)(\xi+a^{-1}m)\neq 0 \quad \text{everywhere on $E$ for all } 0\leq m\leq N-1.$$ Let $I\subseteq \frac{1}{a}\T$ be an interval such that $I+a^{-1}k_0\subseteq E$ for some $k_0\in \Z$. We now prove that the claimed statement holds on $I$ using mathematical induction. 
    Note that since $e^{-2\pi i an_kb^{-\ell}(\xi-\frac{n}{a})}$ is nonzero everywhere for all $k\in \N$ and $\ell\in \N$, we must have $$\text{supp}\bigparen{\Psi_{\frac{1}{a}}(D_{b^\ell}g)}=\text{supp}\bigparen{\Psi_{\frac{1}{a}}(T_{an_i}D_{b^\ell}g)}=\text{supp}\bigparen{\Psi_{\frac{1}{a}}(T_{an_j}D_{b^\ell}g)}$$ for any $1\leq i,j\leq N.$ Assume that $$\text{dim}\bigparen{\clspan\set{\Psi_{\frac{1}{a}}(T_{an_k}D_{b^\ell}g)}}_{n_k\in J_1}\Eq N-1$$
    on $I$ for any $J_1\subseteq \set{n_k}_{k=1}^N$ that contains exactly $N-1$ numbers for some $N\geq 2.$ Without loss of generality, 
    suppose to the contrary that there exists some subset $I_0\subseteq I$ of positive measure at which 
\begin{equation} 
\label{contrapositive_eq_1}
\Psi_{\frac{1}{a}}(T_{an_N}D_{b^\ell}g)\Eq \sum_{k=1}^{N-1}\sigma_k\Psi_{\frac{1}{a}}(T_{an_k}D_{b^\ell}g)
\end{equation}
    for some complex-valued functions $\set{\sigma_k(\xi)}_{k=1}^{N-1}$ defined on $I_0.$ 
    Consequently, we have
    \begin{equation}
    \label{dila_eq_2}
    \Fc\bigparen{T_{an_N}D_{b^{\ell}}g} \Eq \sum_{k=1}^{N-1} \tau_k\Fc\bigparen{T_{an_k}D_{b^{\ell}}g}
    \end{equation}
on $I_0+\frac{1}{a}\Z$ for some $\frac{1}{a}$-periodic functions $\tau_k$ defined on $I_0.$ Note that $\tau_k(\xi)\neq 0$ for almost every $\xi\in E$ and all $1\leq k\leq N-1$. If $\tau_i$ were to vanish on some $I_1\subseteq I_0$ of positive measure, then we would have        $$\Fc\bigparen{T_{an_N}D_{b^{\ell}}g} \Eq \sum_{\substack{1\leq k\leq N-1\\k\neq n_i}}\tau_k\Fc\bigparen{T_{an_k}D_{b^{\ell}}g}\quad \text{a.e. on }I_1+\frac{1}{a}\Z,$$
which contradicts the hypothesis that $\clspan\set{\Psi_{\frac{1}{a}}(T_{an_k}D_{b^\ell}g)}_{n_k\in J_1}$ is $N-1$-dimensional
    on $I$ for any $J_1\subseteq \set{n_k}_{k=1}^N$ that contains exactly $(N-1)$ numbers. Now using the fact that $\tau_k$ are $\frac{1}{a}$-periodic, we obtain that for any $0\leq m\leq N-1.$
    \begin{equation}
    \label{dila_eq_3}
     \Fc\bigparen{T_{an_N}D_{b^{\ell}}g}(\xi+a^{-1}m) \Eq \sum_{k=1}^{N-1}\tau_{k}(\xi) \Fc\bigparen{T_{an_k}D_{b^{\ell}}g}(\xi+a^{-1}m).
     \end{equation}
     By the construction of $I$, $\Fc(D_{b^{\ell}}g)(\xi+a^{-1}m)\neq 0$ on $I_0+a^{-1}k_0$ for all $0\leq m\leq N-1.$
    It follows that 
    \begin{equation}
    \label{dila_eq_4}
    e^{-2\pi i an_Nb^{-\ell}\xi}\Eq \sum_{k=1}^{N-1}\tau_k(\xi)e^{-2\pi i an_kb^{-\ell}\xi}e^{-2\pi i (n_k-n_N)mb^{-\ell} }.
    \end{equation}
    on $I_0+a^{-1}k_0$ for all $0\leq m\leq N-1.$ For notational convenience, we denote $\tau_k(\xi)e^{-2\pi i an_kb^{-\ell}\xi}$ and $e^{2\pi i(n_k-n_N)b^{-\ell}}$ by $f_k(\xi)$ and $a_k$ for all $1\leq k\leq N-1$, respectively. Note that $a_k\notin \set{0,1}$ for all $1\leq k\leq N-1.$ Equation (\ref{dila_eq_4}) now implies that 
    \begin{equation*}
    \sum_{k=1}^{N-1}(a_k^m-1)f_k(\xi)= 0 \quad\text{on }I_0+a^{-1}k_0 \text{ for all $1\leq m\leq (N-1)$.}  
    \end{equation*}
    We then consider the $(N-1)\times (N-1)$ matrix $M$ defined by 
$$M=
\begin{bmatrix}
a_1-1 & \cdots & a_{N-1}-1\\
a_1^2-a_1 & \cdots & a_{N-1}^2-1\\
\vdots & & \vdots\\
a_1^{N-1}-1 & \cdots & a_{N-1}^{N-1}-1
\end{bmatrix}.$$
By dividing $k$-th column by $a_k-1$, we see that the invertibility of $M$ is equivalent to the invertibility of the following matrix        
    $$M'=\begin{bmatrix}
1 & 1 & \cdots & 1\\
1+a_1 & 1+a_2 & \cdots & 1+a_N\\
1+a_1+a_1^2 & 1+a_2+a_2^2 & \cdots & 1+a_N+a_N^2\\
\vdots & \vdots & & \vdots\\
1+a_1+\cdots+a_1^{N-2} &
1+a_2+\cdots+a_2^{N-2} &
\cdots &
1+a_N+\cdots+a_N^{N-2}
\end{bmatrix}.$$
Performing a series of row reductions, we see that $M'$ is indeed a Vandermonde matrix in disguise. As a result, $f_k(\xi)=0$ everywhere on $I_0+a^{-1}k_0$ for all $1\leq k\leq N-1.$ However, this implies that $\tau_k$ vanishes on  $I_0+\frac{1}{a}\Z$ for all $1\leq k\leq N-1,$ which contradicts our hypothesis.
\end{proof}

We are now ready to prove the main result of this subsection.
\begin{theorem}
\label{dila_incomp_semiregu_trans}
    Fix a nonzero real number $b$ with $|b|>1$. Assume that $M\subseteq L^2(\R)$ is a nontrivial closed subspace that is closed under the dilation operator $D_b.$ Then $M$ does not admit a complete set of semi-regular $a$-translates for any $a>0$.
\end{theorem}
\begin{proof} By consider $b^2$ if necessary, we may simply consider the case $b>1.$
    Suppose to the contrary that there exists some $a>0$, $N\in\N$ and $\set{g_i}_{i=1}^N\subseteq L^2(\R)$ such that
    $$M\Eq \clspan{\bigset{S_{a,N}(g_i,J_i,\set{\alpha_{ij}}_{j=1}^{n_i})}}$$
    for some 
    subsets $J_1,\dots,J_N$ of $\Z$ and some finite sets $\set{\alpha_{1j}}_{j=1}^{n_1},\dots,\set{\alpha_{Nj}}_{j=1}^{n_N}$ of $\R$.
    Without loss of generality, we may assume that $g_i\in M$ for all $1\leq i\leq N.$
    By Theorem \ref{dila_fini_linear_inde}, $M$ must be infinite-dimensional. Consequently, without loss of generality, we may assume that $J_1$ is an infinite subset of $\Z.$ 
    Let $L=N+\sum_{i=1}^Nn_i$ and let $J$ be the range function associated with the following $a$-shift-invariant subspace of $L^2(\R)$
    $$V=\clspan{\bigset{S_{a,N}(g_i,\Z)\cup \bigparen{\cup_{\substack{1\leq j\leq n_i\\ 1\leq i\leq N}} S_{a,1}(T_{\alpha_{ij}}g,\Z)}}}.$$ 
    Then clearly, we have that dim$(J(\xi))$ is at most $L$ for almost every $\xi\in \frac{1}{a}\T.$ 
    Since $M$ is closed under dilation, it follows that $\set{T_{ak}D_{b^{\ell}}g_1}_{k\in J_1}\subseteq V$ for all $\ell\in \N.$ However, the closed span of $\set{\Psi_{\frac{1}{a}} (T_{ak}D_{b^{\ell}}g_1)}_{k\in J_1,\,\ell\in\N}$ is infinite-dimensional on some subset of positive measure of $\frac{\T}{a}$ by Lemma \ref{dila_infini_dim_lemma}, which is a contradiction.
\end{proof}
\begin{example}
(a) A \emph{wavelet system} is a sequence in $L^2(\R)$ of the form 
$$\mathcal{W}(\psi,a,b)=\set{a^{n/2}\psi(a^nx-bk)}_{k,n\in \Z},$$
where $\psi \in L^2(\R)$ and $a>1,~b>0$ are fixed. A classic example is the \emph{Haar system} in $L^2(\R)$, which corresponds to the wavelet system generated by the function $\psi=\chi_{[0,\frac{1}{2})}-\chi_{[\frac{1}{2},1)}$ with the dilation parameter $a=2$ and the translation parameter $b=1.$ By Theorem \ref{dila_incomp_semiregu_trans}, it follows that the closed span of any wavelet system does not admit any semi-regulate $a$-translates for any $a>0.$

\smallskip
(b) The \emph{Hardy space} $H^2(\R)=\set{f\in L^2\,|\,\Fc(f)=0\text{ a.e. on }(-\infty,0)}$ is closed under every positive dilation. Consequently, $H^2(\R)$ does not admit any semi-regulate $a$-translates for any $a>0.$

\end{example}
\subsection{Reflection and Fourier Transform} 
We close this section by investigating the compatibility of the existence of semi-regular $a$-translates and the last two operations: reflection and the Fourier transform. 

\medskip
\noindent \textbf{Complete Sets of Semi-Regular $a$-Translates in Reflection-Invariant Subspaces}

We will see below that, not only in the finite-dimensional setting but also in infinite-dimensional setting that there exist closed subspaces that are closed under reflection and simultaneously admit complete sets of semi-regular $a$-translates.

\begin{example}[Finite-dimensional Case] Fix $a>0$ and fix $N\in \N.$
Let $g\in L^2(\R)$ be such that  \begin{equation*}
\begin{aligned}
\widehat{g}(\xi)=\left\{
             \begin{array}{ll}
             1, &\text{if $\xi\in (0,1]$,} \\
             \\
             e^{2\pi ia(N-1)\xi}, &\text{if $\xi\in [-1,0]$,}\\ \\
             0&\text{otherwise.}
             \end{array}
\right.
\end{aligned}
\end{equation*}
Clearly, we have $\widehat{g}=M_{a(N-1)}\Rk\widehat{g}.$
Therefore, we have $T_{ak}g=T_{a(N-1-k)}\Rk g$ for all $0\leq k\leq N-1$. Consequently, the closed subspace $\clspan{\set{T_{ak}g}}_{k=0}^{N-1}$ is closed under reflection. \qeddef
\end{example}


\begin{example}[infinite-dimensional Case]
   Fix $a>0$ and let $g\in L^2(\R)$ be such that \begin{equation*}
\begin{aligned}
\widehat{g}(\xi)=\left\{
             \begin{array}{ll}
             \frac{1}{2}, &\text{if $\xi\in (0,\frac{1}{a}]$,} \\
             \\
            \frac{1}{2} e^{2\pi ia\xi}, &\text{if $\xi\in [-\frac{1}{a},0]$.}\\ \\
0&\text{otherwise.}
             \end{array}
\right.
\end{aligned}
\end{equation*} 
Then we have $T_{ak}\Rk g=T_{-a(k-1)}g$ and $T_{-ak}\Rk g=T_{a(k+1)}g$ for any $k\geq 0.$ So, $\clspan{\set{T_{ak}g}}_{k\in \Z}$ is closed under reflection. In particular, since $\sum_{k\in \Z}|T_{\frac{k}{a}}\widehat{g}|^2=1$ almost everywhere, it follows that $\set{T_{ak}g}_{k\in \Z}$ is an orthogonal basis for its closed span (for example, see \cite[Theorem 10.19]{Hei11}). \qeddef
\end{example}

\noindent\textbf{Complete Sets of Semi-Regular $a$-translates in Fourier-Transform-Invariant Subspaces}

Next, we move on to the case of Fourier transform. Similar to the case of subspaces that are closed under reflection, there exist finite-dimensional subspaces that are closed under Fourier transform. For instance, for each $N\in \N$ the closed span of $N$ linearly independent eigenfunctions of Fourier transform provides such an example. Here an eigenfunction of the 
Fourier transform is a function $f$ in $L^2(\R)$ satisfying $\Fc(f)=\lambda f$ for some $\lambda \in \C.$ Thus, for each $N\geq 1$ and any $a>0$ there exist $N$-dimensional closed subspaces of $L^2(\R)$ that are closed under Fourier transform and also admit complete sets of regular $a$-translates. However, this is not the case when the closed subspace is infinite-dimensional. 
We will show that if $M$ is an infinite-dimensional closed subspace of $L^2(\R)$ that is closed under Fourier transform, then $M$ fails to admit a complete set of semi-regular $a$-translates for any $a>0$ satisfying $a^2\in \Q.$ The case $a^2\notin \Q$ is left open.

\medskip

\noindent \textbf{Solvability of $a$-periodic Equations of Translates II} We have proved that certain $a$-periodic equations of translates are not solvable on $\R$ in Theorem \ref{Non_exist_iterative_eq}. We observe that Theorem \ref{Non_exist_iterative_eq} does not hold if we only consider $a$-periodic equations of translates over a subset of $\R$. For example, let $g=\chi_{[\frac{1}{2},1]}+\chi_{[\frac{3}{2},2]}$. It is then clear that $g=\sum_{k\in \N} T_{k}g$ holds on $(0,\frac{1}{2})+\Z$ as well as on $[\frac{3}{2},2]$. With further assumptions on where the $a$-periodic equation is considered, as well as on the choice of $a$, we extend Theorem \ref{Non_exist_iterative_eq} as follows. 
\begin{theorem}
\label{iterative_eq_lemma_subset}
     Fix $a>0$. Let $g_1,\dots,g_N \in L^2(\R)$ be nonzero functions and let  
    $\set{n_i}_{i=0}^N \subseteq \Z$ be an increasing sequence. Assume that $b>0$ is such that $\frac{a}{b}(n_i-n_0)\in \Z$ for any $1\leq i\leq N.$
    Then for any $b$-periodic functions $\sigma_1,\dots,\sigma_N$, not all identically zero, that are finite almost everywhere the following statements hold.
    \begin{enumerate}
    \setlength\itemsep{0.5em}
        \item [\textup{(a)}] There does not exist a subset $E$ of $\text{supp}(T_{an_0}g)$ with positive measure such that  \begin{equation}
        \label{iterative_equation_subset}
            T_{an_0}g= \sum_{i=1}^N\sigma_i T_{an_i}g \quad \text{a.e. on } E+\frac{1}{b}\Z.
            \end{equation}
        \item [\textup{(b)}] There does not exist a subset $E$ of $\text{supp}(T_{an_N}g)$ with positive measure such that  
        \begin{equation}
        T_{an_N}g= \sum_{i=1}^N\sigma_i T_{an_i}g \quad \text{a.e on } E+\frac{1}{b}\Z. 
        \end{equation}
    \end{enumerate}
    \qeddef
\end{theorem}

We need following lemma to prove this theorem.
\begin{lemma}
\label{generator_lemma_subset}
Fix $a>0$ and let $g$ be a nonzero function in $L^2(\R)$. Assume that $f$ and $h$ are in $\clspan{\set{T_{ak}g}}_{k\in \Z}$. Then we have $\clspan{\set{T_{ak}f}}_{k\in \Z}=\clspan{\set{T_{ak}h}}_{k\in \Z}$ if and only if $\text{supp}(\widehat{f})=\text{supp}(\widehat{h}).$
\end{lemma}
\begin{proof}
($\Rightarrow$) This direction follows from Lemma \ref{generator_lemma}.\medskip

($\Leftarrow$)    Since $f,g\in S(g,\Z)$, by Lemma \ref{nece_condit_complete_regu_trans}, there exist some $\frac{1}{a}$-periodic functions $\tau$ and $\gamma$ that are finite almost everywhere such that 
\begin{equation}
\label{support_eq}
\widehat{f}(\xi)=\tau(\xi)\widehat{g}(\xi)\quad\text{and}\quad \widehat{h}(\xi)=\gamma(\xi)\widehat{g}(\xi)
\end{equation}
We claim that we can make $\text{supp}(\tau)$ equal to $\text{supp}(\gamma)$ by redefining $\tau$ and $\gamma$ if necessary.
 Without loss of generality, suppose that there exists some subset $E\subseteq [0,\frac{1}{a}]$ of positive measure for which $\tau$ vanishes on $E+\frac{1}{a}\Z$, while $\gamma$ is nonzero on the same set.
But $\text{supp}(\widehat{f})$ is equal to $\text{supp}(\widehat{h})$, the only possibility is that $\widehat{g}$ vanishes on $E+\frac{1}{a}\Z$. In this case, we simply redefine $\tau$ by setting $\tau =\gamma$ on $E$. Note that Equation \ref{support_eq} remains valid under this modification. Now, since $\text{supp}(\tau)= \text{supp}(\gamma)$, we can define a $\frac{1}{a}$-periodic function $\sigma(\xi)$ by 

\begin{equation*}
\sigma(\xi) \Eq \begin{cases}
   \dfrac{\gamma(\xi)}{\tau(\xi)}
    & \text{if $\xi\in \text{supp}(\tau)$,} \\[5 \jot]
    \displaystyle 0
    & \text{if $\xi\notin \text{supp}(\tau)$}.
\end{cases}
\end{equation*}
and obtain that $\widehat{h}=\sigma\widehat{f}$. Thus, by Lemma \ref{nece_condit_complete_regu_trans}, we have $h\in S_{a,1}(f,\Z)$. Arguing similarly, we also see that $f\in S_{a,1}(h,\Z).$
\end{proof}

\begin{proof}[Proof of Theorem \ref{iterative_eq_lemma_subset}] (a) Suppose to the contrary that there exist $b$-periodic functions $\sigma_1,\dots,\sigma_N$ that are finite almost everywhere and some subset $E\subseteq \text{supp}(T_{an_0}g)$ of positive measure such that Equation (\ref{iterative_equation_subset}) holds in $E+\frac{1}{b}\Z.$ Applying translation if necessary, we may assume that $E\subseteq [0,b].$ Let $m(x)=\max_{1\leq i\leq N}\sigma_i(x)$ and define the $b$-periodic function $\tau$ on $[0,b]$ by \begin{equation*}
\tau(x) \Eq \begin{cases}
   1
    & \text{if $x\in E\text{ and }m(x)\leq 1$,} \\[5 \jot]
     \frac{1}{m(x)}
    & \text{if $x\in E\text{ and }m(x)\geq 1$,} \\[5 \jot]
    \displaystyle 0
    & \text{if $x\notin E$}.
\end{cases}
\end{equation*}
Therefore, we have 
\begin{equation}
\label{iterative_equation_subset2}
\tau T_{an_0}g=\sum_{i=1}^N\tau\sigma_i T_{an_i}g
\end{equation}
for almost every $x\in \R$. Moreover, we now have $\tau\sigma_iT_{an_i}g\in L^2(\R)$ for all $1\leq i\leq N.$ Applying the inverse Fourier transform to Equation (\ref{iterative_equation_subset2}), we obtain 
\begin{equation}\label{iterative_equation_subset3}
M_{an_0}\Fc^{-1}\bigparen{(T_{-an_0}\tau) g}\Eq \sum_{i=1}^N M_{an_i}\Fc^{-1}\bigparen{(T_{-an_i}\tau) (T_{-an_i}\sigma_i) g}.
\end{equation}
Note that $\Fc^{-1}\bigparen{(T_{-an_i}\tau) g}\in  \clspan\bigset{S_{\frac{1}{b},1}\bigparen{\Fc^{-1}g,\Z}}$ for all $0\leq i\leq N$. Consequently,  for each $1\leq i\leq N$ we have 
\begin{equation*}
\Fc^{-1}\bigparen{(T_{-an_i}\sigma_i)(T_{-an_i}\tau)  g}\in \clspan\bigset{S_{\frac{1}{b},1}\bigparen{\Fc^{-1}\bigparen{(T_{-an_i}\tau) g},\Z}}\subseteq \clspan\bigset{S_{\frac{1}{b},1}\bigparen{\Fc^{-1}g,\Z}}.
\end{equation*}
Because of the assumption $\frac{a}{b}(n_i-n_0)\in \Z$ for all $1\leq i\leq N$, we see that $$\text{supp}\bigparen{(T_{-an_0}\tau) g}\Eq \text{supp}\bigparen{(T_{-a(n_i-n_0)}T_{-an_0}\tau) g}\Eq \text{supp}\bigparen{(T_{-an_i}\tau) g}$$
for all $1\leq i\leq N.$ It then follows that, by Lemma \ref{generator_lemma_subset}, 
\begin{equation}
\label{generator_lemma_subsetIII}
\clspan\bigset{S_{\frac{1}{b},1}\bigparen{\Fc^{-1}\bigparen{(T_{-an_0}\tau) g},\Z}}\Eq \clspan\bigset{S_{\frac{1}{b},1}\bigparen{\Fc^{-1}\bigparen{(T_{-an_i}\tau) g},\Z}}
\end{equation} for all $1\leq i\leq N$. Let $g_0=\Fc^{-1}\bigparen{(T_{-an_0}\tau)g}.$ Combining Equation (\ref{iterative_equation_subset3}) and (\ref{generator_lemma_subsetIII}), we conclude that  
$$M_{an_0}g_0 \in \clspan{\bigset{M_{an_i}T_{\frac kb}g_0\,|\,1\leq i\leq N,k\in\Z}}.$$
Thus, $$M_{an_0}S_{\frac{1}{b},1}(g_0,\Z)\subseteq \clspan{\bigset{M_{an_i}T_{\frac kb}g_0\,|\,1\leq i\leq N,k\in\Z}},$$
which is a contradiction to Lemma \ref{linear_independence_lattice}. Statement (b) can be proved similarly.
\end{proof}
\begin{lemma} \label{contain_finite_F_lemma}
Fix $a,N>0$ with $a^2\in\Q.$
Assume that $M=\clspan{\bigset{S_{a,N}\bigparen{g_i,J_i}}}$ for some $\set{g_i}_{i=1}^N\subseteq L^2(\R)$, not all identically zero, and some subsets $\set{J_i}_{i=1}^N\subseteq \Z$. Then for each $1\leq i\leq N$ there exist at most finitely many $n\in J_i$ for which $\Fc^{-1}(T_{an}g_i)$ is in $M$. 
\end{lemma}

\begin{proof}
We proceed with the proof by contradiction.
 Without loss of generality, we may assume that there is an infinite subset $J_1'$ of $J_1$ for which $\set{\Fc^{-1}(T_{an}g_1)}_{n\in J_1'}\subseteq M.$ Let $J$ be the range function associated with $S_{a,N}(g_i,\Z)$. Then we have $\dim(J(\xi))\leq N$ almost everywhere on  $\frac{1}{a}\T$. By our assumption, we have that $$\clspan{\bigset{\Psi_a\bigparen{\Fc^{-1}(T_{an}g_1}(\xi)}}_{n\in J_1'}\Eq \clspan{\bigset{\bigparen{g_1(\xi-an-k/a)}_{k\in\Z}}}_{n\in J_1'}\subseteq J(\xi)$$  
 for almost every $\xi \in \frac{\T}{a}.$ It remains to show that $$\dim\bigset{\clspan{\bigset{\bigparen{g_1(\xi-an-k/a)}_{k\in\Z}}}_{n\in J_1'}}>N$$ on a subset of $\frac{\T}{a}$ of positive measure to reach a contradiction.
 Since $a^2\in \Q$, there exists an increasing sequence $\set{n_k}_{k\in\N}$ contained in $J_1'$ for which $$(a^2n_i \text{ mod 1})\Eq (a^2n_j\text{ mod 1})\quad \text{for all }i,j\in \N.$$
Consequently, for all $i,j\in \N$ we have $$\text{supp}\bigparen{\Psi_a\bigparen{\Fc^{-1}(T_{an_i}g_1)}}\Eq \text{supp}\bigparen{\Psi_a\bigparen{\Fc^{-1}(T_{an_j}g_1)}}.$$
 Now, suppose to the contrary that there exist some $L\in \N$ with $L\geq 2$ and some subset $E$ of $\text{supp}\bigparen{\Psi_a(\Fc^{-1}(T_{an_1}g_1)}$ of positive measure such that 
 \begin{equation}
 \label{equation_subset_iterative}
\Psi_a\bigparen{\Fc^{-1}(T_{an_L}g_1)}\in \clspan{\bigset{\Psi_a\bigparen{\Fc^{-1}(T_{an_i}g_1)}}}_{i=1}^{L-1}
 \end{equation}
almost everywhere on $E$. Equivalently, there exist some $\frac{1}{a}$-periodic functions $\sigma_1,\dots,\sigma_{L-1}$ for which 
$T_{an_L}g\Eq \sum_{i=1}^{L-1}\sigma_iT_{an_i}g$ on $E+\frac{1}{a}\T.$
However, this is a contradiction to Theorem \ref{iterative_eq_lemma_subset}. Thus, we must have $$\dim\bigset{\clspan{\bigset{\Psi_a\bigparen{\Fc^{-1}(T_{an}g_1)}}}_{n\in J_1'}}=\infty$$
on $\text{supp}\bigparen{\Psi_a\bigparen{\Fc^{-1}(T_{an_1}g_1)}}$, which is again a contradiction. 
\end{proof}

\begin{theorem} 
\label{Fourier_incomp_semiregu_trans}
   Assume that $M$ is an infinite-dimensional closed subspace of $L^2(\R)$ that is closed under Fourier transform. Then $M$ does not admit any complete set of semi-regular $a$-translates for any $a>0$ satisfying $a^2\in \Q.$
\end{theorem}
\begin{proof} Fix $a>0$ with $a^2\in \Q$. 
 Suppose that there exist some $N\in\N$,  some $\set{g_i}_{i=1}^N\subseteq L^2(\R)$, some subsets $J_1,\dots,J_N$ of $\Z$ and some finite subsets $\set{\alpha_{ij}}_{\substack{1\leq j\leq N_i\\1\leq i\leq N}}$ of $\R$ such that
    $$M\Eq \clspan{\set{S_{a,N}(g_i,J_i,\set{\al_{ij}}_{j=1}^{n_i})}}_{i=1}^N$$
    Without loss of generality, we may assume that $J_1$ is infinite. By Lemma \ref{contain_finite_F_lemma}, it remains to consider the case that $\clspan{\set{T_{\alpha_{ij}}g_i}}_{\substack{1\leq j\leq K_i\\1\leq i\leq N}}\not\subseteq \clspan{\bigset{S_a(g_i,I_j)}}_{i=1}^N$ for any subsets $I_1,\dots,I_N$ of $\Z.$ Otherwise, we would have that $M$ fails to be closed under Fourier transform. Since $a^2\in \Q$, we can find some increasing sequence $\set{n_k}_{k\in\N}$ contained in $J_1$ such that $$(a^2n_i \text{ mod 1})\Eq (a^2n_j\text{ mod 1})\text{  for all }i,j\in \N.$$ By the assumption that $M$ is closed under Fourier transform, for each $k\in \N$ there exists some nonzero scalars $\set{c_{ij}^{(k)}}_{\substack{1\leq j\leq N_i\\1\leq i\leq N}}$ and some $\frac{1}{a}$-periodic functions $\set{\sigma^{(k)}_{i}}_{i=1}^N$ that are finite almost everywhere and not all identically zero such that 
        \begin{equation}
        T_{an_k}g_1=\sum_{i=1}^N \sigma^{(k)}_i\widehat{g_i}+\sum_{i=1}^N\sum_{j=1}^{N_i} c_{ij}^{(k)}e_{-\alpha_{ij}}\widehat{g_i}.
        \end{equation}
    It remains to show that $\clspan{\bigset{T_{an_k}g_1-\sum_{i=1}^N \sigma^{(k)}_i\widehat{g_i}}_{k\in \N}}$ is infinite-dimensional to reach a contradiction. Suppose to the contrary that there exists some $L\in \N$ such that for each $\ell \geq L+1$ there exists some nonzero sequence of scalars $(d^{(\ell)}_k)_{k=1}^L$ for which \begin{equation}
        \label{contradic_eq_1}
   T_{an_\ell}g_1-\sum_{i=1}^N\sigma_i^{(\ell)}\widehat{g_i}\Eq \sum_{k=1}^Ld^{(\ell)}_k\bigparen{T_{an_k}g_1-\sum_{i=1}^N \sigma^{(k)}_i\widehat{g_i}}.
    \end{equation}
    Equation (\ref{contradic_eq_1}) then implies $$\bigset{\Fc^{-1}\bigparen{T_{an_\ell}g_1-\sum_{k=1}^L d^{(\ell)}_kT_{an_k}g_1}}_{\ell\geq L+1} \subseteq \clspan{\bigset{S_{a,N}(g_i,\Z)}}.$$ Similarly, we will show that $\clspan{\bigset{\Psi_\frac1a\bigparen{\Fc^{-1}\bigparen{T_{an_\ell}g_1-\sum_{k=1}^L d^{(\ell)}_kT_{an_k}g_1}}}_{\ell\geq L+1}}$ is infinite-dimensional on some subset of $\frac{1}{a}\T$ of positive measure to obtain the desired contradiction. 
    
    We first claim that for all $\ell\geq L+1$ we have 
    \begin{align}
    \begin{split}\text{supp}\bigparen{\Psi_\frac1a\bigparen{\Fc^{-1}\bigparen{T_{an_\ell}g_1-\sum_{k=1}^L d^{(\ell)}_kT_{an_k}g_1}}}\Eq \text{supp}\bigparen{\Psi_\frac1a\bigparen{\Fc^{-1}\bigparen{T_{an_1}g_1}}}.
    \end{split}
    \end{align}
    Since $(a^2n_i \text{ mod 1})\Eq (a^2n_j\text{ mod 1})$, $\textup{supp}(\Psi_\frac1a(\Fc^{-1}(T_{an_i}g_1)))\Eq \textup{supp}(\Psi_\frac1a(\Fc^{-1}(T_{an_j}g_1)))$ for all $i,j\in \N.$
    So, the support of  $\Psi_\frac1a\bigparen{\Fc^{-1}\bigparen{T_{an_\ell}g_1-\sum_{k=1}^Ld^{(\ell)}_kT_{an_k}g_1}}$ is contained in the support of $\Psi_\frac1a\bigparen{\Fc^{-1}(T_{an_1}g_1)}$ for all $\ell\geq L+1.$ If $\Psi_\frac1a\bigparen{\Fc^{-1}\bigparen{T_{an_\ell}g_1-\sum_{k=1}^L d^{(\ell)}_kT_{an_k}g_1}}$ were equal to zero on some subset $E$ of $\Psi_\frac1a(\Fc^{-1}(T_{an_1}g))$ with positive measure for some $\ell\geq L+1$, then we would have $$T_{an_\ell}g_1=\sum_{k=1}^L\tau_kT_{an_k}g_1$$ on $E+\frac{1}{a}\Z$ for some $E\subseteq \frac{1}{a}\T$ of positive measure and some $\frac{1}{a}$-periodic functions $\tau_k$, which is a contradiction to Theorem \ref{iterative_eq_lemma_subset}.
    
    Finally, suppose that there exist some subset $\widetilde{E}$ of $\text{supp}\bigparen{\Psi_\frac1a(\Fc^{-1}(T_{an_{1}}g_1)}$ with positive measure and some integer $K\geq 1$ such that 
\begin{equation}
    T_{an_{L+1+K}}g_1-\sum_{k=1}^L d^{(L+K+1)}_kT_{an_k}g_1=\sum_{\ell=L+1}^{L+K}\eta_\ell \bigparen{T_{an_\ell}g_1-\sum_{k=1}^L d^{(\ell)}_kT_{an_k}g_1}
\end{equation}
almost everywhere on $\widetilde{E}+\frac{1}{a}\Z$ for some $\frac{1}{a}$-periodic functions $\eta_1,\dots,\eta_\ell$ that are finite almost everywhere. Then we would have $$ T_{an_{L+1+K}}g_1=\sum_{k=1}^L d^{(L+K)}_kT_{an_k}g_1+\sum_{\ell=L+1}^{L+K}\eta_\ell \bigparen{T_{an_\ell}g_1-\sum_{k=1}^L d^{(\ell)}_kT_{an_k}g_1}$$
on $\widetilde{E}+\frac{1}{a}\Z$, which is a contradiction to Theorem \ref{iterative_eq_lemma_subset}. We then conclude that $$\text{dim}\bigparen{\clspan{\bigset{\Psi_\frac1a\bigparen{\Fc^{-1}\bigparen{T_{an_\ell}g_1-\sum_{k=1}^L d^{(\ell)}_kT_{an_k}g_1}}}_{\ell\geq L+1}}}=\infty,$$
which is a contradiction.
\end{proof}
\begin{remark}

(a) Arguing similarly to the proof of Theorem \ref{Fourier_incomp_semiregu_trans}, one can also show that if $M$ is an infinite-dimensional closed subspace of $L^2(\R)$ that is closed under Fourier transform, Then $M$ does not admit any complete set of semi-regular $a$-translates with generators whose supports coincide with those of their Fourier transforms, for any choice of $a$.

\medskip
(b)
The results of Section \ref{dilation_section} extend to $\R^d$ for any $d\geq1$, as their proofs of do not rely on the solvability of certain $a$-periodic equations of translates on $\R$. For those results that do rely on such solvability, it remains unclear to us whether they could be generalized to the setting of $L^2(\R^d)$ for any $d>1$. For example, the proof of Theorem \ref{Non_exist_iterative_eq} does not carry over to $\R^d.$  We leave the case of $L^2(\R^d)$ as an open question.

\medskip
(c) The notions of range functions and fiberization maps, which are specific to $L^2(\R)$, play a critical role in Section 3. As a result, our proofs in Section 3 cannot be generalized to $L^p(\R)$ for $p\neq 2$. We also leave the case of $L^p(\R)$ as an open question. 

\medskip
(d) Theorem \ref{modu_incomp_semiregu_trans}, \ref{dila_incomp_semiregu_trans} and \ref{Fourier_incomp_semiregu_trans} are sharp in the sense that 
one cannot extend the number of generators or the number of ``irregular" translates to infinity. For instance, $L^2(\R)$ is closed under modulation, dilation and Fourier transform, but it also admits a complete set consisting of regular translates of infinite many functions (\cite[Theorem 3.3]{Bo00}) and a complete set consisting of infinite irregular translates of a single function (\cite{OL96}). \qeddef
\end{remark}
\section{Operations that are incompatible with frames of translates}
\label{frame_incom_complete}
We investigate in this section the compatibility between frame of translates and the modulation operator in $L^2(\R)$, as well as between frame of translates and the Fourier transform. To maintain the consistency in the underlying space throughout this paper, we will only work primarily in the setting of $L^2(\R)$. However, we note that all results of this section can be easily extended to $L^2(\R^d)$ for any $d>1.$

We begin with a brief introduction to frames. For relatively recent textbook recountings on frame theory, we refer to \cite{Chr16} and \cite{Hei11}. Let $M$ be closed subspace of $L^2(\R)$. A sequence $\set{g_n}\inN\subseteq L^2(\R)$ is said to be a \emph{frame} for $M$ if there exist some positive constants $A\leq B$, called \emph{frame bounds}, such that 
\begin{equation}
\label{frame_equation}
A \norm{f}_{L^2(\R)}^2\Le \sumli |\ip{f}{g_n}|^2\Le B \norm{f}_{L^2(\R)}^2 \quad \text{ for all $f\in M$. }
\end{equation}
Inequality (\ref{frame_equation}) is usually called the \emph{frame inequality}. Since the satisfaction of Inequality (4.1) is independent of the order of a given sequence, we will slightly abuse the terminology by referring to a countable set as a frame if it can be ordered to satisfy the frame inequality.
We say that a sequence $\set{g_n}\inN$ is a \emph{Bessel sequence} in $M$ if at least the upper inequality of Inequality (\ref{frame_equation}) is satisfied. 
It is known that for each frame $\set{g_n}\inN$, there exists another frame $\set{g_n^*}$ in $M$, called \emph{the associated canonical dual frame}, such that 
\begin{equation}
\label{frame_expansions}
f \Eq \sumli \ip{f}{g_n^*}\,g_n,
\end{equation}
with the unconditional convergence of the series in Equation (\ref{frame_expansions}) for all $f\in M.$ Therefore, frames are a more general notion than that of unconditional bases that are bounded above and below in norm, in the sense that every such a unconditional basis is a frame for $L^2(\R)$. In addition to Equation (\ref{frame_expansions}), another fact that follows from the upper frame inequality is the boundedness of the linear operator $R$ from $\ell^2(\N)$ to $L^2(\R)$, called the \emph{associated synthesis operator}, defined by
$R\bigparen{(c_n)\inN}=\sumli c_ng_n.$
If $\set{g_n}\inN$ is a Bessel sequence in $L^2(\R)$, then it is necessary that $\norm{R}\leq B^{1/2}$, where $B$ is an upper frame bound of $\set{g_n}\inN.$ 

We are particularly interested in two specific types of frames in this section: \emph{frames of translates} and \emph{Gabor frames} (or \emph{Weyl-Heisenberg frames}). A frame of translates is a frame that consists of at most 
countably many translates of finitely many functions in $L^2(\R)$.  
A \emph{Gabor system} is a subset of $L^2(\R)$ that is of the form 
$$\Gc_N(g_i,\Gamma_i)=\bigset{M_wT_xg_i\,|\,(x,w)\in \Gamma_i,1\leq i\leq N},$$
where $\set{g_i}_{i=1}^n\in L^2(\R)$ are the associated generators and $\set{\Gamma_i}_{i=1}^N\subseteq \R^2$ are subsets that contains all associated time-frequency shifts. 
If $\Gc_N(g_i,\Gamma_i)$ is a frame for some closed subspace $M$ of $L^2(\R)$, then we refer to it as a \emph{Gabor frame} for $M$ generated by $g_1,\dots,g_N$ to emphasize its structure. 
Consequently, frames of translates are just special cases of Gabor frames in which no frequency shifts are involved. To be consistent in the use of notations throughout this section, we will write $\Gc_N(g_i,\Gamma_i\times \set{0})$ to denote the system of translates generating by $g_1,\dots,g_N$ and the sets $\Gamma_1,\dots,\Gamma_N\subseteq \R.$ Another term that will be used frequently in this section is the \emph{basic sequences} (resp.\ \emph{unconditional basic sequences}), which means sequences that are Schauder bases (resp.\ unconditional bases) for their closed spans.



Next, we provide a brief overview of \emph{modulation spaces}. For more details and related applications of modulation spaces, we refer to \cite{Gro01} and \cite{BO20}. \begin{definition}\label{Modulation_definition}
Fix a nonzero Schwartz function $\psi\in S(\R)$. 
\begin{enumerate}\setlength\itemsep{0.5em}
    \item [\textup{(a)}] The \emph{short-time Fourier transform} of $f\in L^2(\R)$, denoted by $V_{\psi}f$, is the complex-valued measurable function on $\R^2$ defined by  $$V_{\psi}f(x,w)\Eq \overline{\ip{M_wT_x \psi}{f}} \Eq  \ip{f}{M_wT_x \psi}.$$
     \item [\textup{(b)}] For $1\leq p \leq 2$, the \emph{modulation space} $M^{p}(\R)$ is the space consisting of 
      $f\in L^2(\R)$ for which $\norm{V_\psi f}_{L^{p}(\R^2)}$ is finite, i.e.,
     $$M^{p}(\R)=\bigset{f \in L^2(\R)\,\big|\,  \norm{f}_{M^{p}(\R)}\Eq \norm{V_\psi f}_{L^{p}(\R^2)}<\infty}.$$
      \qeddef
    \end{enumerate}
\end{definition} 
Modulation spaces were introduced by Feichtinger in the early 1980s and were further developed in a series of collaborations with Gr\"{o}chenig.  Since then, modulations spaces have been recognized as the appropriate spaces of the study of time-frequency analysis. Loosely speaking, modulation spaces are a family of function spaces which classify functions based on their joint time and frequency localization (via short-time Fourier transform).  
It is worth mentioning that, given the main focus of this paper, the definition of modulation spaces that is presented in this paper is merely a specific family of modulation spaces that are commonly found in the literature.

We will see in the Theorem \ref{modulation_plus_M1_no_ft} that the deciding factor in whether the property of being closed under modulation can coexist with a frame of translates in a closed subspace of $L^2(\R)$ is whether the closed subspace contains a function with good joint time and frequency localization. To prove this result, we need following two lemmas. 
\begin{lemma}
\label{unif_discrete_sum_finite}
    Fix $1\leq p\leq 2$ and fix $1\leq q\leq \frac{2p}{3p-2}$.
    Let $f\in M^p(\R)$ and let $g\in M^q(\R)$ be two functions. Then for any $d>0$  we have $$\sum_{n,k\in\Z}\bignorm{V_fg\cdot\chi_{Q_{d,n,k}}}_{L^\infty(\R^2)}^{\frac{pq}{p+q-pq}}<\infty.$$
    Here $Q_{d,n,k}$ denotes $[dn,d(n+1)]\times [dk,d(k+1)].$
   
\end{lemma}
\begin{proof}
    Since both $f$ and $g$ are in $M^2(\R)=L^2(\R)$, the short-time Fourier transform $V_gf(x,w)$ is continuous on $\R^2$, and hence $V_gf(x,w)$ is well-defined for every $(x,w)\in \R^2.$ Fix $\gamma\in S(\R)$ with $\ip{\gamma}{\gamma}=1.$ 
    We first show that $|V_gf(x,w)|\leq \bigparen{|V_\gamma f|\ast|V_g \gamma|}(x,w)$ for all $(x,w)\in \R^2.$
    Let $\set{g_n}\inN$ and $\set{f_n}\inN$ be two sequences in $S(\R)$ that converge to $g$ and $f$ in $M^2(\R)$, respectively. Then by \cite[Lemma 11.3.3]{Gro01} we have $$|V_{g_n}f_n|(x,w)\Le \bigparen{|V_{\gamma}f_n|\ast |V_{g_n}\gamma|}(x,w).$$
    Using the Triangle Inequality and Cauchy–Bunyakovsky–Schwarz inequality, we obtain 
    \begin{align*}
    \begin{split}
        \norm{V_{g_n}f_n-V_gf}_\infty&\Le \norm{V_{g_n}(f_n-f)}_\infty+\norm{V_{g_n-g}f}_\infty\\
        &\Le \norm{g_n}_{M^2(\R)}\norm{f_n-f}_{M^2(\R)}+\norm{g_n}_{L^2(\R)}\norm{f_n-f}_{M^2(\R)}.
        \end{split}
    \end{align*}
So, $V_{g_n}f_n(x,w)\rightarrow V_gf(x,w)$ as $n$ tends to $\infty$ for every $(x,w)\in \R^2.$ On the other hand, using the Triangle Inequality and Young's convolution inequality, we see that 
\begin{align*}
    \begin{split}
       \bignorm{|V_{\gamma}f_n|\ast |V_{g_n}\gamma|-|V_{\gamma}f|\ast |V_{g}\gamma|}_\infty &\Le \bignorm{\bigabs{V_{\gamma}(f_n-f)}\ast |V_{g_n}\gamma|}_{L^\infty}+\bignorm{|{V_{\gamma}f}|\ast |V_{g_n-g}\gamma|}_{L^\infty}\\
       &\Le \norm{f_n-f}_{M^2(\R)}\norm{g_n}_{M^2(\R)}+\norm{f}_{M^2(\R)}\norm{g_n-g}_{M^2(\R)}.
    \end{split}
\end{align*}
Therefore, $\bigparen{|V_{\gamma}f_n|\ast |V_{g_n}\gamma|}(x,w)$ converges to $\bigparen{|V_{\gamma}f|\ast |V_{g}\gamma|}(x,w)$ for every $(x,w)\in \R$ as $n$ tends to $\infty.$ Thus, we have $|V_gf(x,w)|\leq \bigparen{|V_\gamma f|\ast|V_g \gamma|}(x,w)$ for all $(x,w)\in \R^2.$

Next, let $C>0$ be a constant depending only on $\delta$, $p$ and $q$ such that $$\sum_{n,k\in \Z}\norm{V_gf\cdot\chi_{Q_{d,n,k}}}_\infty^{^{\frac{pq}{p+q-pq}}}\Le C \sum_{n,k\in \Z}\norm{V_gf\cdot\chi_{Q_{1,n,k}}}_\infty^{^{\frac{pq}{p+q-pq}}}.$$
Then by \cite[Theorem 11.1.5]{Gro01} and \cite[Theorem 12.2.1]{Gro01}, we see that \begin{align*}
    \begin{split}
        \sum_{n,k\in \Z}\norm{V_gf\cdot_{Q_{1,n,k}}}_\infty^{^{\frac{pq}{p+q-pq}}}
        &\Le  C \sum_{n,k\in \Z}\bignorm{\bigparen{|V_\gamma f|\ast |V_g \gamma|\ast |V_{\gamma}\gamma|}\cdot\chi_{Q_{1,n,k}}}^{\frac{pq}{p+q-pq}}_\infty \\
        &\Le C \Bigparen{\bigparen{\sum_{n,k\in \Z}\bignorm{|V_{\gamma}\gamma|\cdot\chi_{Q_{1,n,k}}}_\infty} \norm{f}_{M^p(\R)} \norm{g}_{M^q(\R)}}^{\frac{pq}{p+q-pq}}.\qedhere
    \end{split}
\end{align*}
 \end{proof}

We remark that the double projection method that will be used in the proof of following Lemma has been used in \cite{RS95}, \cite{CDH99} and \cite{DH00}. For notational convenience, we denote by $B_R(x,w)$ the ball in $\R^2$ centered at $(x,w)$ with radius $R.$ For any $y,z\in \R^2$, we use $\norm{y-z}$ to denote the Euclidean distance between $y$ and $z.$ Finally, for a subset $\Gamma\subseteq \R^d$ we write $|\Gamma|$ to denote its cardinality.

\begin{lemma} \label{homogeneous_distance_lemma} Fix $1\leq p\leq 2.$ 
Let $M$ be an infinite-dimensional closed subspace of $L^2(\R)$ that contains an unconditional basic sequence of the form $\Gc_N(\phi_i,\Lambda_i)$ for some $\phi_1,\dots,\phi_N \in M^p(\R)$ and some subsets $\Lambda_1,\dots,\Lambda_N \subseteq\R^2.$ 
 Assume that $M$ admits a Gabor frame $\Gc_K(g_i,\Gamma_i)$ with generators $g_i\in M^{q_i}(\R)$, where $1\leq q_i\leq \frac{2p}{3p-2}$ for all $1\leq i\leq K.$
 Then for each $\epsilon>0$ there exists some $R_0>0$ such that for any $(x,w)\in \bigcup\limits_{i=1}^N\Lambda_i$, any $h>0$ and any $R\geq R_0$ we have
\begin{equation}
    \label{cardinality_comparison_equation}
(1-\epsilon)\sum_{i=1}^N|B_h(x,w)\cap\Lambda_i|\Le\sum_{i=1}^K|B_{R+h}(x,w)\cap\Gamma_i|.
\end{equation}
\end{lemma}
\begin{proof}
    Given that $\Gc_N(\phi_i,\Lambda_i)$ and $\Gc_K(g_i,\Gamma_i)$ are Bessel sequences, each $\Lambda_i$ and each $\Gamma_i$ can be decomposed into finitely many subsequences that are uniformly discrete by \cite[Theorem 3.1]{CDH99}. By decomposing $\Lambda_i$ and $\Gamma_i$ into subsequences if necessary, we may assume that each $\Lambda_i$ and each $\Gamma_i$ are uniformly separated and let $d=\min\limits_{1\leq i\leq N}\inf\set{|x-y|\,|\,x,y\in \Gamma_i}.$ Since $\Lambda_i$ and $\Gamma_i$ are uniformly discrete, for any $(x,w)\in \cup_{i=1}^N\Lambda_i$ and any $R>0$, the following two subspaces 
    $$S_R(x,w)=\clspan{\set{g^*_{i,(a,b)}\,|\,(a,b)\in \Gamma_i\cap B_R(x,w) \text{ for some }1\leq i\leq K}},$$
    and 
     $$W_R(x,w)=\clspan{\set{M_bT_a\phi_i\,|\,(a,b)\in \Lambda_i\cap B_R(x,w) \text{ for some }1\leq i\leq N}},$$
     are finite-dimensional. Here $\set{g^*_{i,(a,b)}}_{\substack{(a,b)\in \Lambda_i\\1\leq i\leq K}}$ denotes the canonical dual frame associated with $\Gc_K(g_i,\Gamma_i).$  
     We then observe that for any $h,R>0$ and any $(x,w)\in \Lambda$ the right-hand side of Equation (\ref{cardinality_comparison_equation}) dominates the rank of the orthogonal projection operator $P_{S_{R+h}(x,w)}$ from $L^2(\R)$ onto $S_{R+h}(x,w)$.  Moreover, since $\Gc_N(\phi_i,\Lambda_i)$ is an unconditional basic sequence, we have the following inequality that
     \begin{equation}
         \label{trace_inequal}
     |B_{h}(x,w)\cap \Lambda|\Le \sum_{i=1}^N\sum_{(a,b)\in B_{h}(x,w)\cap \Lambda_i}|\ip{M_{b}T_a\phi_i}{\phi_{i,a,b}^*}|.
     \end{equation}
     where $\set{\phi_{i,a,b}^*}_{\substack{(a,b)\in \Lambda_i\\1\leq i\leq N}}$ denotes the coefficient functionals associated with $\Gc_N(\phi_i,\Lambda_i).$

    Next, to relate the rank of $P_{S_{R+h}(x,w)}$ to the quantity on the right-hand side of Inequality (\ref{trace_inequal}), we consider the finite-rank linear operator $T$ from  $W_{h}(x,w)$ to $W_{h}(x,w)
$ defined by $T=P_{W_{h}(x,w)}P_{S_{R+h}(x,w)}.$ Clearly, we have $$\textup{trace}(T)\Eq \text{sum of all eigenvalues of $T$}\Le \text{Rank}(T).$$
On the other hand, we compute
\begin{align*}
\begin{split}
\text{trace}(T)&\Eq \sum_{i=1}^N\sum_{(a,b)\in B_h(x,w)\cap \Lambda_i} \ip{T(M_bT_a\phi_i)}{P_{W}\phi^*_{i,a,b}}
\\
&\Eq \sum_{i=1}^N\sum_{(a,b)\in B_h(x,w)\cap \Lambda_i} \bigip{P_{S}(M_bT_a\phi_i)}{P_W\phi^*_{i,a,b}}\\
&\Eq \sum_{i=1}^N \sum_{(a,b)\in B_h(x,w)\cap \Lambda_i} \Bigparen{\bigip{(P_{S}-I)(M_bT_a\phi_i)}{P_{W}\phi^*_{i,a,b}}+ \bigip{M_bT_a\phi}{P_{W}\phi^*_{i,a,b}}}\\
&\Le \sum_{i=1}^N\sum_{(a,b)\in B_h(x,w)\cap \Lambda_i} \bigparen{1-(\sup_{\substack{a,b\in\Lambda_i\\1\leq i\leq N }}\norm{P_W\phi^*_{i,a,b}}_{L^2(\R)})\norm{(P_S-I)M_bT_a\phi_i}_{L^2(\R)}}
\end{split}
\end{align*}
Here, for notational convenience, we denote $P_{W_h(x,w)}$ and $P_{S_{R+h}(x,w)}$ by $P_W$ and $P_S$, respectively. Note that $\sup\limits_{\substack{a,b\in\Lambda_i\\1\leq i\leq N }}\norm{P_W\phi^*_{i,a,b}}_{L^2(\R)}=C<\infty$ since $\Gc_N(\phi_i,\Lambda_i)$ is an unconditional basic sequence. Therefore, it remains to show that for each $\epsilon>0$ there exists some $R_0>0$ such that 
$$\norm{(P_{S_{R+h}(x,w)}-I)M_bT_a\phi_i}_{L^2(\R)}<\frac{\epsilon}{C}.$$
for any $R\geq R_0$, any $h>0$, any $(x,w)\in \cup_{i=1}^N\Lambda_i$, any $(a,b)\in B_h(x,w)\cap(\cup_{i=1}^N\Lambda_i)$ and all $1\leq i\leq N.$

    Let $\epsilon>0$ be an arbitrary real number. If $(a,b)\in B_h(x,w)$ for some $h>0$, then we have $S_R(a,b)\subseteq S_{R+h}(x,w)$. Since $\set{g_{j,a,b}}_{\substack{(a,b)\in \Gamma_j\\
    1\leq j\leq K}}$ is a frame for $M$ , for each $(a,b)\in B_h(x,w)$ and $1\leq i\leq N$ we have 
    \begin{equation}
    \label{stop_equation}
        M_bT_a \phi_i\Eq \sum_{j=1}^K\sum_{(\gamma,w)\in \Gamma_j} \bigip{M_bT_a\phi}{M_{w}T_{\gamma}g_j}\,g_{j,\gamma,w}^*.
    \end{equation}
    It then follows that   
\begin{align*}
\begin{split}
\label{distance_estimate}
\norm{(P_{S_{R+h}(x,w)}-I)M_bT_a\phi_i}_{L^2(\R)}&\Le  \norm{(P_{S_R(a,b)}-I)M_bT_a\phi_i}_{L^2(\R)}
\\
&\Le \Bignorm{\sum_{j=1}^K\sum_{(\gamma,w)\in \Gamma_j\setminus B_R(a,b)} \bigip{M_bT_a\phi_i}{M_{w}T_{\gamma}g_j}\,g_{j,\gamma,w}^*}_{L^2(\R)}\\
&\Le B\, \sum_{j=1}^K\sum_{(\gamma,w)\in \Gamma_j\setminus B_R(a,b)} \Bigabs{\bigip{M_bT_a\phi_i}{M_{w}T_{\gamma}g_j}}^2,
\end{split}
\end{align*}
where $B$ is an upper frame bound of $\set{g_{j,a,b}}_{\substack{(a,b)\in \Gamma_j\\
    1\leq j\leq K}}$. Note that, by definition, we have $|\ip{M_bT_a\phi_i}{M_{w}T_{\gamma}g_j}|=|V_{\phi_i}g_j(a-\gamma,w-b)|.$ Thus, we have 
\begin{align*}
\begin{split}
B\sum_{j=1}^K\sum_{(\gamma,w)\in \Gamma_j\setminus B_R(a,b)} \Bigabs{\bigip{M_bT_a\phi_i}{M_{w}T_{\gamma}g_j}}^2
&\Eq B\, \sum_{j=1}^K\sum_{(\gamma,w)\in \Gamma_i\setminus B_R(a,b)} \bigabs{V_{\phi_i}g_j(a-\gamma,w-b)}^2\\
&\Le B\sum_{j=1}^K\sum_{\substack{n,k\in \Z\\ Q_{d,n,k} \nsubseteq B_{R-d}(0)}} \bignorm{V_{\phi_i}g_j\cdot \chi_{Q_{d,n,k}}}_{L^\infty(\R)}^2
\end{split}
\end{align*}
where $Q_{d,n,k}=[nd,(n+1)d)\times [kd,(k+1)d) $ 
By Lemma \ref{unif_discrete_sum_finite}, there exists some $R_0>0$ such that for all $R\geq R_0$, the final term in the estimate above is less than $\frac{\epsilon}{BC}$ for all $1\leq i\leq N$. 
 
\end{proof}

We have an immediate corollary of Lemma \ref{homogeneous_distance_lemma} as follows.
For any $x\in \R$ and any subset $\Gamma\subseteq \R$ the distance from $x$ to $\Gamma$, denoted by $d(x,\Gamma)$ is defined by $$d(x,\Gamma)=\inf_{y\in \Gamma}\norm{x-y}.$$ 

\begin{corollary}
\label{distance_criterion}
    Fix $1\leq p\leq 2.$
Let $M$ be an infinite-dimensional closed subspace of $L^2(\R)$ that contains an unconditional basic sequence of the form $\Gc_N(\phi_i,\Lambda_i)$ for some $\phi_1,\dots,\phi_N $ in $M^p(\R)$ and some subsets $\Lambda_1,\dots,\Lambda_N \subseteq\R^{2}.$ 
 Assume that $M$ admits a Gabor frame $\Gc_K(g_i,\Gamma_i)$ with $g_i\in M^{q_i}(\R)$ for some $1\leq q_i\leq \frac{2p}{3p-2}$ for all $1\leq i\leq N.$
Then 
$\bigparen{\sup\limits_{x\in \Lambda_j} d(x,\cup_{i=1}^K\Gamma_i)}<\infty$.
for all $1\leq j\leq N.$
\end{corollary}


We are now ready to prove the main theorem of this section.
\begin{theorem}
\label{modulation_plus_M1_no_ft}
    Let $M\subseteq L^2(\R)$ be a closed subspace. Assume that $M$ contains a closed subspace $M'$ that is closed under modulation and $M'\cap M^1(\R)\neq \set{0}$. 
    Then $M$ does not admit a frame of translates.
\end{theorem}
\begin{proof}
Let $b$ be a nonzero real number such that $M_b(M')\subseteq M'$.
Suppose to the contrary that $M$ admits a frame of the form $\Gc_N(g_i,\Gamma_i\times \set{0})$ for some $g_1,\dots,g_N\in L^2(\R)$ and some subsets $\Gamma_1,\dots,\Gamma_N\subseteq \R.$

By Corollary \ref{distance_criterion}, it suffices to find an unconditional basic sequence contained in $M'$ of the form $\Gc_1(\phi,\Lambda)$ for some subset $\Lambda\subseteq \set{0}\times \R$.
Let  $\phi\in M\cap M^1(\R)$ be a nonzero function. We first choose $k_0$ large enough that \begin{equation}
    \label{Wiener_amalgam_frame_condition}
A\Le \sum_{k\in \Z} \bigabs{\phi(\cdot-(bk_0)^{-1} k)}^2\Le B.
\end{equation}
 almost everywhere on $\R$ for some $A,B>0$. By \cite[Corollary 3.5]{Wa92} there exists some $a_0>0$ such that if $a<a_0$, then $\Gc_1(\phi,(bk_0)^{-1},a)$ is a Gabor frame for $L^2(\R).$ Specifically, one can choose $a_0>0$ to be some number that satisfies $$\sum_{\substack{n\in \Z \\n\neq 0}}\Bignorm{\sum_{\substack{k\in \Z}}\phi(\cdot-(bk_0)^{-1}k)\overline{\phi}(\cdot-na_0-(bk_0)^{-1} k)}_{L^\infty(\R)}<A.$$
almost everywhere on $\R$. By Ron-Shen Duality Principle (\cite[Theorem 2.2]{RS97}), $\Gc_1(\phi, a^{-1},bk_0)$ is an unconditional basic sequence. Consequently, $M'$ contains an unconditional basic sequence of the form $\Gc_1(\phi, 0,bk_0)$,
 which is a contradiction.
\end{proof}

Theorem \ref{modulation_plus_M1_no_ft} is a significant improvement over one of the main results in \cite{CDH99}. To illustrate this, we present the following example.
\begin{example}
    Fix $a\geq 0$ and $b>0.$ Let $\phi=e^{-\frac{\pi x^2}{2}}$ be the Gaussian function. Then every closed subspace of $L^2(\R)$ containing $\Gc_1(\phi,a,b)$ does not admit a frame of translates. We first consider the case that $a\neq 0.$ It is known that $\Gc_1(\phi,a,b)$ is a frame for $L^2(\R)$ if $ab<1$; $\Gc_1(\phi,a,b)$ is a complete set in $L^2(\R)$ if $ab=1$; $\Gc_1(\phi,a,b)$ is an unconditional basic sequence in $L^2(\R)$ if $ab>1$. Consequently, the case $ab\leq 1$ follows directly from \cite[Theorem 1.2]{CDH99}, while the case $ab>1$ is now implied by Corollary \ref{distance_criterion}. Finally, the case $a=0$ follows from the fact that $\Gc_1(\phi,0,b)$ is an unconditional basic sequence for any $b>0$. \qeddef
\end{example}

\begin{remark} 

(a) There is another proof of Theorem \ref{modulation_plus_M1_no_ft} via \emph{Feichtinger Conjecture} (now known as \emph{Feichtinger Theorem}). 
It was conjectured by Feichtinger that every frame that is bounded below in norm can be decomposed into finitely many unconditional basic sequences. Later, Casazza et al.\ proved that Feichtinger Conjecture is equivalent to the statement that every Bessel sequence that is bounded below in norm can be partitioned into finitely many unconditional basic sequences \cite{CCLV05}. Following that, Casazza et al.\ demonstrated that Feichtinger's Conjecture is equivalent to the Kadison-Singer Conjecture \cite{CT06}. Kadison-Singer Conjecture was ultimately resolved by Marcus, Spielman and Srivastava in \cite{MSS15}. For an expository account of Kadison-Singer Conjecture, we refer to \cite{Bo18}.

We now outline how to prove Theorem \ref{modulation_plus_M1_no_ft} by using Feichtinger Theorem. 
Let $Q=[0,1].$ It is known that $M^1(\R)\subseteq W(L^\infty,\ell^1)$, where $$W(L^\infty,\ell^1)\Eq \Bigset{f\colon\R\rightarrow \C \text{ measurable}\,\Big|\,\sum_{k\in \Z}\norm{f\cdot T_{k}\chi_Q}_{L^\infty(\R)}<\infty\,}.$$
Using Walnut Representation(see \cite{Wa92} and \cite[Section 11.5]{Hei11}), one can see that for any fixed nonzero $\phi\in M^1(\R)$ and any fixed $a,b>0$ the corresponding Gabor system $\Gc_1(\phi,a,b)$ is a Bessel sequence. 
    Building upon this fact, we can now apply Feichtinger Theorem to decompose $\Gc_1(\phi,a,b)$ into finitely many unconditional basic sequences, from which we extract the desired unconditional basic sequence that consists solely of the frequency shifts of $\phi$.  However, this approaches loses information about the specific structure of the desired unconditional basic sequence contained in $M$, such as which frequency shifts are involved. For research papers regarding the structure of the decomposition yielded by Feichtinger Theorem, we refer to \cite{BCMS19} and \cite{Bo24}.
    \medskip

 (b)  
 Assume that $M$ is a closed subspace of $L^2(\R)$ that is invariant under modulation, i.e., there exists some nonzero $b\in \R$ such that $M_b(M)=M.$ Then $\Fc(M)$ is $b$-shift-invariant. It was shown by Bownik in \cite[Theorem 3.3]{Bo00} that every $1$-shift-invariant subspace $S$ admits a countable set of generators.
 That is, there exists some countable subset $J$ of $\N$ and a countable subset $\set{\phi_n}_{n\in J }$ of $L^2(\R)$ such that
 $S=\clspan{\set{T_k\phi_n\,|\,k\in \Z,\,n\in J}}$.
 Using the standard dilation technique (see Proposition \ref{necessary_cond_range}), one can extend this result to $b$-shift-invariant subspaces. Now let $\set{\phi_n}_{n\in J}$ be a set of generators for $\Fc(M)$. Using Proposition \ref{necessary_cond_range}, it follows that $M$ does not admit any frame of translates if there exists some function $\phi$ in $M^1(\R)$ such that $$\Psi_{ \frac{1}{b}}(\phi)(\xi)\in\clspan{\set{\Psi_{ \frac{1}{b}}(\phi_n)(\xi)\,|\,n\in J}}$$
 for almost every $\xi\in \frac{1}{b}\T.$ This idea also applies to closed subspaces that are invariant under both modulation and translation, where a different type of range functions come into play. For the characterization of modulation-translation-invariant subspaces using range functions, we refer the reader to \cite{Bo07}.
 
\end{remark}

The next corollary shows that if an infinite-dimensional closed subspace of $L^2(\R)$ is closed under Fourier transform and simultaneously admit a frame of translates, then not all generators of that frame of translates can be well-localized in both time and frequency.
\begin{corollary}
Let $M$ be an infinite-dimensional closed subspace of $L^2(\R)$ that is closed under Fourier transform. Assume that $M$ admits a frame of translates $\Gc_N(g_i,\Gamma_i\times\set{0})$ for some $g_1,\dots,g_N\in L^2(\R)$ and some subsets $\Gamma_1,\dots,\Gamma_N\subseteq\R$. Then 
     there exists some $1\leq i_0\leq N$ such that $g_{i_0}\notin M^p(\R)$ for any $1\leq p\leq \frac{4}{3}.$

\end{corollary}
\begin{proof}
    By assumption, $\Gc_{N}\bigparen{\widehat{g_i},\set{0}\times (-\Gamma_i)}$ is also a frame for $\Fc(M)$. By Feichtinger Theorem, for each $1\leq i\leq N$ there exists an infinite subset $\Gamma_i'$ of $\Gamma_i$ such that $\Gc_{1}\bigparen{\widehat{g_i},\set{0}\times (-\Gamma_i') }$ is an unconditional basic sequence.     
    Next, suppose to the contrary that $g_i\in M^p(\R)$ for some $1\leq p\leq \frac{4}{3}$ and all $1\leq i\leq N.$ Since modulation spaces are invariant under Fourier transform, we see that $\widehat{g_i}\in M^p(\R)$ for some $1\leq p\leq \frac{4}{3}$ and all $1\leq i\leq N.$ Statement (a) then follows from  Corollary \ref{distance_criterion} and the fact that $1\leq p\leq \frac{2p}{3p-2}$ if and only if $1\leq p\leq \frac{4}{3}.$
\end{proof}

We conclude this paper by investigating the compatibility between modulation, Fourier transform and the existence of Schauder basis of translates. Recall that if $\set{g_n}\inN$ is a Schauder basis for $L^2(\R)$ with associated coefficient functionals $\set{g_n^*}\inN$, then there exists some constant $C>0$, called the associated \emph{basis constant}, such that $$1\Le \norm{g_n}_{L^2(\R)}\norm{g_n^*}_{L^2(\R)}\Le C$$  
(see, for example, \cite[Section 4.3]{Hei11}). Consequently, if $\set{g_n}\inN$ is a Schauder basis of translates for $L^2(\R)$, then it is necessary that there exist some constants $A,B>0$ such that 
$$A\Le \sup_{n\in\N}\norm{g_n^*}\Le B.$$
    \begin{theorem} \label{compati_modulation_Schauder_basis} 
Let $M$ be a closed subspace of $L^2(\R)$. Assume that $M$ admits a Schauder basis of translates $\Gc_N(g_i,\Gamma_i\times\set{0})$ for some $g_1,\dots,g_N\subseteq L^2(\R)$ and some $\Gamma_1,\dots,\Gamma_N$.
Assume that one of the following condition holds. 
\begin{enumerate}
\setlength\itemsep{0.5em}
    \item [\textup{(a)}] There exists some closed subspace $M'\subseteq M$ that is closed under modulation and $M'\cap M^1(\R)\neq \set{0}.$
    \item [\textup{(b)}] $M$ is infinite-dimensional and is closed under Fourier transform. 
\end{enumerate}
 Then there exists some $1\leq i_0\leq N$ such that $g_{i_0}\notin M^1(\R).$
\end{theorem}
\begin{proof} 
Suppose to the contrary that $g_i\in M^1(\R)$ for all $1\leq i\leq N$.

\medskip
(a)
Let $\phi$ be a nonzero function that belongs to $M'\cap M^1(\R)$. Arguing similarly to the proof of Theorem \ref{modulation_plus_M1_no_ft}, we see that there exists some infinite uniformly discrete subset $\Lambda\subseteq \R$ such that $\Gc_1(\phi,\set{0}\times \Lambda)$ is an unconditional basic sequence. To make the proof of statement (a) applicable to statement (b), we may assume that $\Gc_1(\phi,\set{0}\times \Lambda)$ is a basic sequence. 

By \cite[Lemma 4.1]{DH00}, each of $\Gamma_i$ is uniformly discrete. By passing to subsequences if necessary, we may assume that each $\Gamma_i$ is uniformly discrete and let $$d=\min\limits_{1\leq i\leq N}\inf\set{|x-y|\,|\,x,y\in \Gamma_i}.$$ 
Similarly, given any $w\in \Lambda$ and any $R>0$, we define the following two finite-dimensional subspaces 
    $$S_R(0,w)=\clspan{\set{g^*_{i,\gamma}\,|\,(\gamma,0)\in \Gamma_i\cap B_R(0,w) \text{ for some }1\leq i\leq N}},$$
    and 
     $$W_R(0,w)=\clspan{\set{M_b\phi\,|\,(0,b)\in \Lambda\cap B_R(0,w) }},$$
    where $\set{g^*_{i,\gamma}}_{\substack{(\gamma,0)\in \Lambda_i\\1\leq i\leq N}}$ denotes the associated coefficient functionals of $\Gc_N(g_i,\Gamma_i\times\set{0}).$  
Let $0<\epsilon<1$ and repeat the proof of Lemma \ref{homogeneous_distance_lemma} up to Equation (\ref{stop_equation}). That is, \begin{equation}
\label{stop_equation_II}
M_b \phi\Eq \sum_{i=1}^N\sum_{(\gamma,0)\in \Gamma} \bigip{M_b\phi}{T_{\gamma}g_i}\,g_{i,\gamma}^*,\quad \text{for all $b\in \Lambda$.}
\end{equation}
  It then follows that for any $R,h>0$, any $w \in \Lambda$ and any $(0,b)\in B_h(0,w)\cap (\set{0}\times \Lambda)$ and all $1\leq i\leq N,$
\begin{align}
\begin{split}
\label{repeat_estimate}
\norm{(P_{S_{R+h}(0,w)}-I)M_b\phi_i}_{L^2(\R)}&\Le  \norm{(P_{S_R(0,b)}-I)M_bT_a\phi_i}_{L^2(\R)}
\\
&\Le B\, \sum_{j=1}^N\sum_{(\gamma,0)\in \Gamma_j\setminus B_R(0,b)} \bigabs{\bigip{M_b\phi}{T_{\gamma}g_j}}\\
&\Eq B\, \sum_{j=1}^N\sum_{(\gamma,0)\in \Gamma_j\setminus B_R(0,b)} \bigabs{V_{\phi}g_j(-\gamma,b)},\\
\end{split}
\end{align}
where $B=\sup \set{\norm{g_{i,\gamma}^*}_{L^2(\R)}\,|\,\gamma \in \Gamma_i, 1\leq i\leq N}$. Since $$\sum_{j=1}^N\sum_{(\gamma,0)\in \Gamma_i\setminus B_R(0,b)} \bigabs{V_{\phi}g_j(-\gamma,b)}
\Le \sum_{j=1}^N\sum_{\substack{n,k\in \Z\\ Q_{d,n,k} \nsubseteq B_{R-d}(0)}} \bignorm{V_{\phi}g_j\cdot \chi_{Q_{d,n,k}}}_{L^\infty(\R)},$$
it follows that, by Lemma \ref{unif_discrete_sum_finite}, there exists some $R_0$ such that the last term in the estimate above is sufficiently small whenever $R\geq R_0.$ Here  $Q_{d,n,k}$ denotes $[nd,(n+1)d)\times [kd,(k+1)d).$
Thus, we have proved that given any $0<\epsilon<1$ there exists some $R_0>0$ such that for all $R\geq R_0$, $h>0$ any $w\in \Lambda$ the following statement holds: 
\begin{equation}
   \label{cardinality_ineqaulity} 
(1-\epsilon)|B_h(0,w)\cap\Lambda|\Le\sum_{i=1}^N|B_{h+R}(0,w)\cap\Gamma_i|.
\end{equation}
Since $\min\limits_{1\leq i\leq N}\sup\set{d(x,\Gamma_i)\,|\,x\in \Lambda}=\infty$, we obtain a contradiction.

(b) Assume that condition (b) holds. Then being closed under Fourier transform implies that $M$ simultaneously admits a Schauder basis of translates and a Schauder basis that consists solely of frequency shifts of $\widehat{g_i}$, which is a contradiction to Equation (\ref{cardinality_ineqaulity}).
\end{proof}

\section*{acknowledgement}
We thank Ryszard Szwarc for bringing Theorem \ref{dila_fini_linear_inde} to our attention. We are also grateful to Marcin Bownik for his valuable discussions on this project. Finally, we thank the anonymous comments that helped us identify and correct serious errors in an earlier version of this manuscript.


\end{document}